\documentclass[11pt,letterpaper,final]{amsart}

\usepackage[pdfauthor={D. Schrittesser, A. Tornquist},
            pdftitle={The Ramsey property and higher dimensional mad families},
            pdfproducer={Latex with hyperref}]{hyperref}

\usepackage{todonotes}

\usepackage[utf8]{inputenc}
\usepackage[T1]{fontenc}
\usepackage{lmodern}
\usepackage[english]{babel}
\usepackage{csquotes}

\usepackage{amsmath}
\usepackage{amstext}
\usepackage{amssymb}
\usepackage{amsthm}
\usepackage[mathb,matha]{mathabx}

\usepackage{mathrsfs}
\usepackage[all]{xy}
\usepackage{datetime}

\usepackage[notref,notcite]{showkeys}

\usepackage{datetime}

\usepackage[shortlabels]{enumitem}
\usepackage{accents}
\newlength{\dhatheight}

\makeatletter
\newcommand{\whathat}[1]{{%
  \mathpalette\double@widehat{#1}%
}}
\newcommand{\double@widehat}[2]{%
  \sbox\z@{$\m@th#1\widehat{#2}$}%
  \ht\z@=.9\ht\z@
  \widehat{\box\z@}%
}
\makeatother

\theoremstyle{plain}
\newtheorem{thm}{Theorem}[section]
\newtheorem{lem}[thm]{Lemma}

\newtheorem{cor}[thm]{Corollary}
\newtheorem{claim}[thm]{Claim}

\theoremstyle{definition}
\newtheorem{dfn}[thm]{Definition}

\newtheorem{fact}[thm]{Fact}
\newtheorem{facts}[thm]{Facts}
\newtheorem{rem}[thm]{Remark}

\theoremstyle{remark}

{\end{minipage}\end{equation}}

\newenvironment{eqpar*}{\begin{equation*}\begin{minipage}{0.8\columnwidth}}%
{\end{minipage}\end{equation*}}

\newcommand{\AC}{\ensuremath{\mathsf{AC}}}

\mathchardef\mhyphen="2D
\newcommand{\RUnif}{\mathsf{R\mhyphen Unif}}

\newcommand{\tD}{\mathbf}

\DeclareMathOperator{\splitset}{sp}
\DeclareMathOperator{\splitenum}{\widehat{sp}}

\newcommand{\consec}[1][A]{\mathbin{\prec^{#1}}}

\DeclareMathOperator{\powerset}{\mathcal{P}}

\newcommand{\fin}[1][]{{\mathsf{Fin}^{#1}}}
\newcommand{\finp}[1][]{{\mathsf{Fin}^{#1 +}}}
\newcommand{\finpp}[1][]{{\mathsf{Fin}^{#1 ++}}}
\newcommand{\finalph}{\fin[\laddermap]}
\newcommand{\finalphp}{{\finalph}^+}
\newcommand{\finalphpp}{{\finalph}^{++}}
\newcommand{\finsq}{\fin[2]}
\newcommand{\finsqp}{{\fin[2]}^{+}}
\newcommand{\finsqpp}{{\fin[2]}^{++}}

\newcommand{\subsetfin}[1][\laddermap]{\mathbin{\subseteq_{\fin[#1]}}}

\newcommand{\from}[1][z]{\widetilde{#1}}

\newcommand{\fromx}[2][{}]{\widetilde{#2}^{#1}}

\newcommand{\laddermap}{\gamma}
\newcommand{\ladderres}[1][s]{\laddermap | {#1}}

\DeclareMathOperator{\higherdom}{\overline{dom}}

\providecommand{\reals}{\mathbb{R}}
\providecommand{\nat}{\mathbb{N}}

\DeclareMathOperator{\dom}{dom}
\DeclareMathOperator{\lh}{lh}

\DeclareMathOperator{\ran}{ran}

\providecommand{\res}{\mathbin{\upharpoonright} }
\providecommand{\conc}{ \mathbin{{}^\frown}}

\providecommand{\ZF}{\ensuremath{\mathsf{ZF}}}
\providecommand{\HOD}{\mathbf{HOD}}

\newcommand{\eL}{L}

\providecommand{\CCR}{\ensuremath{\mathsf{CC}_{\reals}}}
\providecommand{\DC}{\ensuremath{\mathsf{DC}}}
\providecommand{\DCR}{\ensuremath{\mathsf{DC}_\reals}}
\providecommand{\AD}{\ensuremath{\mathsf{AD}}}
\providecommand{\PD}{\ensuremath{\mathsf{PD}}}

\providecommand{\setdef}{\;|\;}

\def\N{{\mathbb N}}

\newcommand{\cond}{\ensuremath{\Gamma}}
\newcommand{\condargs}{\ensuremath{\Gamma(\tD d, t'_0, t'_1, s, m)}}

\author[Schrittesser]{David Schrittesser}
\address{David Schrittesser, Institute for Advanced Study in Mathematics, Harbin Institute of Technology, 92 West Da Zhi Street, Harbin City, Hei\-long\-jiang Province 150001, China, \emph{and}}

\address{Harbin Institute of Technology Suzhou Research Institute, Building K, 500 Nan Guandu Road, Suzhou City, Jiangsu Province 215104, China.}
\email{david@logic.univie.ac.at}

\author[Törnquist]{Asger Törnquist}

\address{Asger Törnquist, Department of Mathematical Sciences, University of Copenhagen, Universitetspark 5, 2100 Copenhagen, Denmark}
\email{asgert@math.ku.dk}

\title{The Ramsey property and higher dimensional mad families}

\subjclass[2020]{03E05, %
03E15 %
03E60 %
05D10
}

\keywords{Higher mad families, Borel ideals, Ramsey regularity, Fubini-product, Axioms of Determinacy, Solovay's model}

\begin{document}

\begin{abstract}
We prove that under a principle of Ramsey regularity there are no infinite maximal almost disjoint families with respect to the transfinitely iterated Fréchet ideals. The results of the present paper were announced by the authors in the Proceedings of the National Academy of Sciences of the U.S.A.
\end{abstract}

\maketitle

\section{Introduction}

The present paper is a contribution to the area of infinite combinatorics related to maximal almost disjoint (``mad'') families, which has seen rapid growth in recent years \cite{chan-jackson, horowitz-shelah-no-mad, horowitz-shelah-vector, horowitz-todorcevic, neeman-norwood, toernquist}.  
The goal here is to develop more fully the connection between mad families and a combinatorial notion of measurability, viz., the Ramsey Property; and to supply proofs to the theorems that were announced in \cite{pnas}. Specifically, we show that the results obtained in \cite{pnas} for the ideal $\fin$ of finite subsets of $\nat$ (also known as the Fréchet ideal) can be extended to the class of ideals on $\nat$ obtained from $\fin$ by iterating a type of Fubini-product operation (see Definition~\ref{d.iterated.Fréchet} and its preceding paragraph).  
Such ideals were to our knowledge first studied with motivations in harmonic analysis by Kahane \cite{kahane}, and taken up again by Jackson-Kechris-Louveau \cite{jackson-kechris-louveau} in their study of Fréchet-amenability. We study them in the context of almost-disjointness with the intention to enrich the evolving structure theory of ideals on $\nat$ (a small sample of recent contributions to this area are \cite{dobrinen,chodounsky,brendle} where the reader will also find further references).
Since for $k \in \nat$, a $k$-fold Fubini product of $\fin$ is (in a natural way) an  ideal on $\N^k$, we speak of \emph{higher dimensional mad families}.

\medskip

\subsection*{(A) Background} 
If $\mathcal I$ is an ideal on $\N$ and $A,B\subseteq\N$, then we will say that $A$ and $B$ are \emph{$\mathcal I$-almost disjoint} if $A\cap B\in\mathcal I$. Thus, the notion of being $\fin$-almost disjoint 
is the classical notion of almost disjointness. A family $\mathcal A\subseteq\mathcal P(\N)\setminus\mathcal I$ will be called an \emph{$\mathcal I$-almost disjoint family} if $A\cap B\in\mathcal I$ for all $A,B\in\mathcal A$ for which $A\neq B$. A \emph{maximal} $\mathcal I$-almost disjoint family (called an $\mathcal I$-mad family, for short) is an $\mathcal I$-almost disjoint family which is maximal under inclusion. \emph{Finite} $\mathcal I$-mad families exist for trivial reasons (e.g., $\mathcal A=\{\N\}$), whereas the existence of \emph{infinite} $\mathcal I$-mad families follows by a routine application of Zorn's lemma, and so apparently entails the use of the Axiom of Choice ($\AC$).

Let us state Mathias' classical theorem, which started the subject we are investigating:

\begin{thm}[A.R.D.\ Mathias, \cite{mathias-thesis, mathias}]{\ }

(A) There are no analytic infinite $\fin$-mad families.

(B) If $\ZF\,+$``there is a Mahlo cardinal'' is consistent, then $\ZF\,+$``there are no infinite $\fin$-mad families'' is consistent.
\end{thm}

The previous theorem grew out of Mathias' investigation of infinitary Ramsey theory, \emph{happy families}, and Mathias forcing. On this background, Mathias asked two key questions: (1) \emph{Are there any infinite $\fin$-mad families in Solovay's model?} (2) \emph{If all sets of reals have the Ramsey Property, does it follow that there are no infinite $\fin$-mad families?} Since Mathias had proved that all sets have the Ramsey property in Solovay's model, a positive answer to (2) would imply a negative answer to (1), but both questions resisted solution until Question (1) was solved in 2014 by the second author in \cite{toernquist}. Shortly thereafter, Neeman and Norwood \cite{neeman-norwood} gave a different solution to Question (1). Then Horowitz and Shelah showed in \cite{horowitz-shelah-no-mad} that no large cardinal assumption is needed to obtain a model of $\ZF\,+$``there are no infinite $\fin$-mad families''. 

Regarding Question (2), in 2019, a solution was given by the authors in \cite{pnas}:

\begin{thm}[$\ZF\, +\DCR +\RUnif$; \cite{pnas}]\label{t.main1}
If all sets have the Ramsey property, then there are no infinite $\fin$-mad families.
\end{thm}

Here, $\ZF\, +\DCR\, +\RUnif$ denotes the background theory in which the theorem is proved: $\ZF$ denotes denotes Zermelo-Fraenkel set theory \emph{without} $\AC$, while $\DCR$ and $\RUnif$ denote the weak choice principles \emph{Dependent Choice for subsets of $\reals$} and \emph{Uniformization on Ramsey positive sets}, respectively; see Definition~\ref{d.Runif} below.

\subsection*{(B) Higher dimensional mad families} The purpose of the present paper is to prove a version of Theorem \ref{t.main1} for $\mathcal I$-mad families, where $\mathcal I \in \mathfrak F$, the smallest class containing $\fin$ which is closed under taking Fubini products:

\medskip

Recall that if $\mathcal I$ is an ideal on $X$ and $\mathcal J$ is an ideal on $Y$, then the simple \emph{Fubini product} $\mathcal I\otimes\mathcal J$ of $\mathcal I$ and $\mathcal J$ is the ideal on $X\times Y$ defined by
$$
\mathcal I\otimes\mathcal J=\{A\subseteq \N\times \N: \{x\in Y: A(x)\notin \mathcal J\}\in\mathcal I\},
$$
where $A(x)=\{y\in Y: (x,y)\in A\}$. The \emph{finitely} iterated Fréchet ideals are then $\fin^2=\fin\otimes\fin$, $\fin^3=\fin\otimes\fin^2$, and so on. 

For the \emph{transfinitely} iterated Fréchet ideals (cf.\ \cite{kahane,jackson-kechris-louveau}) one needs a more general type of product:
Given an ideal $\mathcal I$ on a set $X$ and ideals $\mathcal J_x$ on sets $Y_x$, where $x$ ranges over $X$, the ideal $\mathcal I\otimes(\mathcal J_x)_{x\in X}$ on the set $\bigcup_{x\in X} \{x\}\times Y_x$ is given by
$$
\mathcal I\otimes(\mathcal J_x)_{x\in X}=\{A\subseteq \bigcup_{x\in X} \{x\}\times Y_x:\{x\in X: A_x\notin\mathcal J_x\}\in\mathcal I\}.
$$
We will call $\mathcal I\otimes(\mathcal J_x)_{x\in X}$ the \emph{Fubini product over $\mathcal I$ of} the family $(\mathcal J_x)_{x\in X}$.

\begin{dfn}\label{d.iterated.Fréchet}
Denote by $\mathfrak F$ the smallest class of ideals such that $\fin \in \mathfrak F$ and $\mathfrak F$ is closed under Fubini products over $\fin$.
\end{dfn}

While the ideals in $\mathfrak F$ could be thought of as ideals on $\N$, in \S\ref{s.alphaD}, we give a useful construction of these ideals where the domain is taken to be the set of leaves of a tree. The set of ideals $\mathfrak F$ is large in some respects; for instance, the Borel hierarchy complexity of the ideals in $\mathfrak F$ is unbounded in rank.

We can now state our main theorem:

\begin{thm}[$\ZF\, +\DCR +\RUnif$]\label{t.main2}
Let $\mathcal I\in\mathfrak F$. If all sets have the Ramsey property, then there are no infinite $\mathcal I$-mad families.
\end{thm}

The proof of Theorem \ref{t.main2} follows a strategy similar to the proof of Theorem \ref{t.main1} in \cite{pnas} in that we reduce the theorem to the case of analytic $\mathcal I$-mad families. We then associate to a given infinite $\mathcal I$-almost disjoint family $\mathcal A$ a function (the \emph{``tilde operator''}). This function sends each infinite subset $W$ of $\N$ to an $\mathcal I$-positive set $\from[W]$ while satisfying certain invariance and pigeon-hole principles. We then define a family of ``search trees'' of attempts at finding an element of $\mathcal A$ which intersects $\from[W]$ in an $\mathcal I$-positive set. Under the assumption that all sets are completely Ramsey, the invariance properties of these search trees, combined with the invariance and pigeon-hole properties of the tilde operator, allow us to find a Ramsey co-null set of $W$ such that $\{\from[W]\} \cup \mathcal A$ is $\mathcal I$-almost disjoint, proving that $\mathcal A$ is not maximal. Though this strategy is similar to the one employed to prove Theorem \ref{t.main1} in \cite{pnas} (the $\fin$ case), the proof is more difficult for the higher-dimensional iterates of $\fin$ because the invariance properties we obtain are, in a sense that later will be made precise, only conditional.

\subsection*{(C) Localization to pointclasses}

Theorem \ref{t.main2} localizes to pointclasses
 in the sense of \cite{moschovakis}. This localization allows us to obtain as a corollary many of the previously known results about higher dimensional mad families that were proved under various regularity assumptions. \emph{For simplicity, in this paper we always assume that pointclasses $\Gamma$ contain all Borel sets, are closed under finite unions and intersections, and continuous preimages}.

\begin{thm}[\ZF]\label{t.main3}
Let $\Gamma$ be a pointclass and assume $\RUnif^*(\Gamma)$. Then there are no Dedekind infinite $\mathcal I$-mad families in $\Gamma$, for any ideal $\mathcal I \in \mathfrak F$.
\end{thm}

\begin{cor}[\ZF]\label{c.intro} 
Let $\mathcal I\in\mathfrak F$. Then the following hold:
\begin{enumerate}
\item There are no infinite analytic $\mathcal I$-mad families.
\item There are no infinite $\mathcal I$-mad families in Solovay's model. In particular, it holds that if $\ZF + $``there is an inaccessible cardinal'' is consistent, so is the theory $\ZF + \DC +$``there are no infinite $\mathcal I$-mad families, for any ideal $\mathcal I \in \mathfrak F$''.
\item Assuming Projective Determinacy $(\PD)$ and $\DCR$, there are no infinite projective $\mathcal I$-mad families.
\item\label{AD} Assuming the Axiom of Determinacy $(\AD)$, there are no infinite $\mathcal I$-mad families in $\eL(\reals)$.
\item\label{ADplus} Assuming $\AD^{+}$, a natural strengthening of the Axiom of Determinacy, there are no infinite $\mathcal I$-mad families.
\end{enumerate}
\end{cor}

Results as in this corollary have recently been generalized to mad families on uncountable cardinals $\kappa$ by Chan and Jackson
\cite{chan-jackson}.

Items (1) and (2) of the corollary for $\mathcal I = \fin$ were first proved in \cite{toernquist}, and a form of (3), as well as (4) and (5) were proved in \cite{neeman-norwood}. 
The results for $\mathcal I\in\mathfrak F$ in the previous corollary were proved using forcing methods in \cite{haga-schrittesser-toernquist}. Note that \eqref{AD} follows from \eqref{ADplus} via Hugh Woodin's result that $\AD\; + V = L(\reals)$ implies $\AD^{+}$.
Whether $\AD$ and $\AD^{+}$ are, in fact, equivalent, is a famous and longstanding unresolved question asked by Woodin. For $\mathcal I=\fin$, the inaccessible cardinal assumption in (2) in the corollary can be taken away, by a result of Horowitz and Shelah \cite{horowitz-shelah-no-mad}. In all likelihood, the previous corollary is not optimal; indeed, we conjecture the inaccessible cardinal can be taken away here too, but at present we do not know how to prove this.

\medskip

\subsection*{Organization of the paper}
To give the reader a sense of orientation, we first prove our main theorem for the case $\mathcal I=\fin^2$.
This is carried out in~\S\ref{s.2D}.
After that, we tackle the general case in~\S\ref{s.alphaD}. 
This latter section is almost completely self-contained; 
a brave reader can skip~\S\ref{s.2D} if desired, but will be asked to refer to~\S\ref{s.2D} 
for a few proofs where the general case differs only in notation from the case $\mathcal I=\fin^2$.

\subsection*{Acknowledgments} The first author thanks the 
FWF for support through Ben Miller's project P29999 and Vera Fischer's START Prize Y1012,
as well as the Government of Canada’s New Frontiers in Research Fund (NFRF), for support through grant NFRFE-2018-02164. Asger Törnquist received funding from the Independent Research Fund Denmark through the grants ``Automorphisms and invariants of operator algebras'' and ``Operator algebras, groups, and quantum spaces''.
We thank Natasha Dobrinen, Stevo Todorčević, Michael Hrušák, Juris Steprāns, and the participants of the Delta Workshop in Wuhan, 2023, for interesting conversations on topics related to this paper.

\section{Background material and notation}\label{s.prelim}

This section serves as a repository for notation and basic definitions and conventions. Unless where it is otherwise explicitly stated, the background theory for the entire paper is $\ZF$ set theory \emph{without} the Axiom of Choice. The principle of Dependent Choices $\DCR$, which is much weaker than the full Axiom of Choice, will be used for some theorems; it is clearly indicated in the text when $\DCR$ is assumed.

\subsection{Almost disjoint families for ideals.} Let $\mathcal I$ be an ideal, in the set-theoretical sense, on a countable set $S$: %
That is, $\mathcal I\subseteq\powerset(S)$ such that for all $I,J\in\powerset(\mathcal \nat)$ it holds that 
(i) $I,J \in \mathcal I \Rightarrow I\cup J \in \mathcal I$ and 
(ii) if $I\subseteq J$ and $J\in \mathcal I$, then also $I\in\mathcal I$.
We write $\mathcal{I}^+$ for $\powerset(S)\setminus \mathcal I$, the \emph{co-ideal} corresponding to $\mathcal I$.
We write $A \subseteq_{\mathcal I} B$ to mean $A\setminus B \in \mathcal I$.

An \emph{$\mathcal I$-almost disjoint family} is a set $\mathcal A \subseteq \mathcal{I}^+$ 
such that any two distinct $A,B \in \mathcal A$ are $\mathcal I$-almost disjoint, 
i.e., $A\cap B \in \mathcal I$.
Such a family is called \emph{maximal} if it is not a proper subset of an almost disjoint family. 
We abbreviate ``maximal $\mathcal I$-almost disjoint family'' by ``$\mathcal I$-mad family''.

If we take $S=\nat$ and $\mathcal I$ to be $\fin$, i.e., the ideal of finite subsets of $\nat$, 
we speak simply of an almost disjoint, resp.\ mad, family.
The reader will notice that we don't follow the usual convention of making ``infinite'' 
a part of the definition of mad family.

\subsection{Descriptive set theory and trees}

For basic concepts of descriptive set theory, such as Polish spaces, Borel sets, etc., we refer to \cite{kechris}. We use the notion of a \emph{pointclass} in the way it is used in \cite{moschovakis}, that is, a pointclass $\Gamma$ is a class of subsets of Polish spaces, e.g. the class of Borel sets in Polish spaces. We always assume that the pointclasses we consider are closed under continuous preimages, finite unions and intersections, and contain the Borel sets as a subclass.

 For $A$ a set, we write $A^{\nat}$ for the set of functions from $\N$ to $A$, and $A^{<\N}$ for the set of finite functions $s:n\to A$, where $n=\{0,1\ldots, n-1\}$.
We follow the descriptive set theoretical conventions regarding trees  (compare \cite[2.C]{kechris}):
By a \emph{subtree of $\nat^{<\nat}$} we mean a set
 $
 \tD c \subseteq \nat^{<\nat}
 $
 which is a tree in the usual sense, that is,
 for every $s \in \tD c$ and $k < \lh(s)$, also $s \res k \in \tD c$.
If $\tD X$ is a subtree of $\nat^{<\nat}$,
 then by a subtree of $\tD X$ we mean a subset of $\tD X$ which is a subtree of $\nat^{<\nat}$.
 We follow the usual abuse of terminology whereby a tree $T$ on $2 \times A$, where $A$ is a set, is a subset of $2^{<\nat}\times A^{<\nat}$
such that $(u, v) \in T \Rightarrow \lh(u) = \lh(v)$ and $(u \res k, v \res k ) \in T$ for $k < \lh(u)$.
When $T$ is a tree on $2\times A$ and $t,t' \in T$, we write $t \subseteq t'$ to mean $u \subseteq u' \land v \subseteq v'$,
where $t = (u,v)$ and $t' = (u',v')$; i.e.,
$t'$ \emph{extends} $t$ in the tree order.
Given $t \in T$, 
\[
T_t = \{t' \in T \setdef  t\subseteq t'\lor t' \subseteq t\}.
\]
As usual we write $[T]$ for the set of branches through $T$,
\[
[T] = \{(x_0, x_1) \setdef (\forall n\in\nat)\; (x_0 \res n, x_1 \res n) \in T\}
\]
and
\[
\pi[T] = \{x_0 \in 2^\nat  \setdef (\exists x_1) \; (x_0, x_1) \in [T]\}.
\]
Throughout, we shall not distinguish between $2^\nat $ and $\powerset(\nat^{<\nat} )$ and treat elements of 
$\pi[T]$ as if they were subsets of $\nat^{<\nat} $ (just fix your favourite bijection $\nat \cong \nat^{<\nat} $ and identify sets with
their characteristic functions).

In order to state precisely our further theorems, we now give formal definitions of some key concepts. ; and the reader who is looking for a quick introduction to the Ramsey theory of $[\nat]^\infty$ should consult \cite[19.D]{kechris}.

\begin{dfn}~\label{d.ramsey}
\begin{enumerate}[(a)]
\item A set $X\subseteq [\N]^\infty$ is said to be \emph{Ramsey} (or have the \emph{Ramsey Property}) if there is $A\in [\nat]^\infty$ such that either $[A]^\infty\subseteq X$ or $[A]^\infty\cap X=\emptyset$. 

\item Given $A\in [\nat]^\infty$ and $a\in [\nat]^{,\infty}$ with $\max(a) < \min(A)$, write as usual $[a,A] = \{S \in [\nat]^\infty \setdef a \subseteq S \subseteq a \cup A\}$ (this is a basic open neighborhood in the Ellentuck topology).

\item A set $X\subseteq [\N]^\infty$ is said to be \emph{completely Ramsey} if for every $A\in [\nat]^\infty$ and $a\in [\nat]^{<\infty}$ with $\max(a) < \min(A)$ there is $B\in [A]^\infty$ such that either $[a,B]^\infty\subseteq X$ or $[a,B]^\infty\cap X=\emptyset$. 

\item A function $f:[A]^\infty\to Y$, where $A\in [\N]^\infty$ and $Y$ is a Polish space, is \emph{Ramsey measurable} (resp.,\ \emph{completely Ramsey measurable}) if the preimage of any open subset of $Y$ is Ramsey (resp., completely Ramsey).
\end{enumerate}
\end{dfn}

\begin{dfn}[Principle of Ramsey positive uniformization, $\RUnif$]\label{d.Runif}\ 
\begin{enumerate}
\item In general, a relation $R\subseteq X\times Y$ is said to have \emph{full projection} if for every $x\in X$ there is $y\in Y$ such that $(x,y)\in R$.
\item $\RUnif$ is the statement: For every relation $R\subseteq [\nat]^\infty \times [\nat]^\infty$ with full projection, there is $A\in [\nat]^\infty$ and a function $f:[A]^\infty\to [\nat]^\infty$ which \emph{uniformizes} $R$ on $[A]^\infty$, that is
$$
(\forall x\in [A]^\infty)\  (x,f(x))\in R.
$$
\item For $\Gamma$ a pointclass, $\RUnif(\Gamma)$ is the statement: For every $\Gamma$-relation $R\subseteq [\nat]^\infty \times [\nat]^\infty$ with full projection, there is $A\in [\nat]^\infty$ and a function $f:[A]^\to [\nat]^\infty$ which uniformizes $R$ on $[A]^\infty$.
\item For $\Gamma$ a pointclass, $\RUnif^*(\Gamma)$ is the same statement as (3), but with the additional requirement that $f$ be completely Ramsey measurable.
\end{enumerate}
\end{dfn}

Given $A\subseteq \nat$, we write $[A]^\infty$ for the set of infinite subsets of $A$ and $[A]^{<\infty}$ for the set of finite subsets of $A$. For $A\subseteq\nat$, write $\hat A: |A|\to\N$ for the strictly increasing function with $\ran(\hat A)=A$.

\section{The two-dimensional case}\label{s.2D}

In this section, we take the first step toward proving our main result: We prove Theorems \ref{t.main2} and \ref{t.main3} in the 2-dimensional case, i.e. where $\mathcal I=\fin^2$. That is, we will prove

\begin{thm}[\DCR]\label{t.2D.precise}
Let $\Gamma$ be a pointclass and assume $\RUnif^*(\Gamma)$. Then there are no infinite $\finsq$-mad families in $\Gamma$.
\end{thm}

We have singled out the case $\mathcal I=\fin^2$ for special treatment, since the proof is easier to follow in this case, yet retains many of the key features of the fully general theorems. The 1-dimensional case $\mathcal I=\fin$ of Theorem \ref{t.2D.precise} was proved in \cite{pnas}, but many issues that make the higher dimensional cases harder don't appear in the case $\mathcal I=\fin$.

\medskip

We point out the following immediate corollary of Theorem \ref{t.2D.precise}, which was proved using forcing methods in \cite{haga-schrittesser-toernquist}:
\begin{cor}\label{c.2D.no.analytic}
There are no analytic $\finsq$-\emph{mad} families.
\end{cor}

Theorem~\ref{t.2D.precise} will be a consequence of a more technical result, Theorem~\ref{t.2D.analytic}, stated further below. Theorem~\ref{t.2D.analytic} is a statement about a certain operator, called \emph{the tilde operator}, which will be defined in subsection \ref{ss.tilde2} below. Given an infinite $\finsq$-almost disjoint family $\mathcal A$, the tilde operator transforms elements of $\finp$ to elements of the co-ideal $(\mathcal A \cup \finsq)^+$ in a way that has many nice properties, and it plays a crucial role in the proofs, in dimension 2 as well as in higher dimensions.

\medskip

{\it Notation\footnote{The notation introduced now, especially $S=\N^2$, is intended to foreshadow the notation used for the fully general case in the next section.}:} Let $S = \nat^2$.
For $\tD X \subseteq S$, we write\footnote{We use boldface letters for subsets of $S$ and lightface letters for subsets of $\nat$.}
\begin{align*}
\tD X(n) &= \{m\in\nat \setdef (n,m) \in \tD X\},\\
\dom(\tD X) &=  \{n\in\nat \setdef (\exists m \in \nat)\; (n,m) \in \tD X\}.
\end{align*}
Moreover, let
\begin{align*}
\dom_\infty(\tD X) &=  \{n\in\nat \setdef \tD X(n) \in \finp\},\\
\tD X^{++} &= \tD X \cap \big(\dom_\infty(\tD X) \times \nat\big),\\
\finsqpp &= \big\{\tD X \in \finsqp \setdef (\forall m) \; \big[\tD X(m)\neq \emptyset \Rightarrow  \tD X(m)\in\finp\big]\big\}.
\end{align*}
With this notation,
$
\fin^2=\{\tD X \subseteq S\mid\dom_\infty(\tD X)\in\fin\}$,
and 
$$
\tD X \in \finsqp \iff \tD X^{++} \in \finsqpp \iff \dom_\infty(\tD X) \in \finp.
$$ 
Moreover, $\tD X \in \finsqpp$ if and only if $\dom_\infty(\tD X) = \dom(\tD X)$ and $\dom(\tD X) \in \finp$.

\medskip

{\it The notation $\hat{\tD{X}}^l(m,n)$ for enumerating $\tD{X}\in\finsqpp$ is defined as follows:} 
Recall from Section \ref{s.prelim} that when $X\subseteq\N$, we denote by $\hat{X}$  the strictly increasing enumeration $\hat{X}\colon \lvert X \rvert \to X$. We extend this notation to $\tD{X}\in\finsqpp$ by defining $\hat{\tD{X}}^l(m,\cdot)$ to be the strictly increasing enumeration of the $m$'th vertical of $\tD{X}$; that is, if $m' = \widehat{\dom(\tD{X})}(m)$, then
\[
\hat{\tD{X}}^l(m,n)=(m',\widehat{\tD{X}(m')}(n)).
\]

\subsection{The two-dimensional tilde operator}\label{ss.tilde2}
For the purpose of this section, fix a $\finsq$-a.d.\  family $\mathcal A$.
The tilde operator will be defined relative to an appropriate sequence $\langle \tD{Z}^l \setdef  l\in\nat\rangle$, chosen to
satisfy the following:
\begin{enumerate}[label=\textup{(\Alph*$_2$)}]
\item \label{Z.refine} For each $l\in \nat$, $\tD{Z}^l \in [\tD A]^{++}_{\finsq}$ for some $\tD A \in \mathcal A$.
\item \label{Z.disjoint} For each $m \in \nat$ there is at most one $l \in \nat$ such that $\tD{Z}^l (m) \neq\emptyset$.
\end{enumerate}
Note that by \ref{Z.disjoint} we also have $\tD{Z}^l \cap \tD{Z}^{l'} =\emptyset$ for $l \neq l'$.

Let us verify that such a sequence is indeed available to us.
\begin{lem}\label{l.Z.2D}
There is a sequence $\langle \tD{Z}^l \setdef  l\in\nat\rangle$ satisfying \ref{Z.refine} and \ref{Z.disjoint} above.
\end{lem}
\begin{proof}
Using $\DCR$ choose a sequence $\langle \tD{A}^l\mid l\in\N\rangle$ of elements of $\mathcal A$ (for the current lemma, $\CCR$ suffices)
and let $B_l = \dom_\infty(\tD{A}^l)$.

There is a sequence $\langle B'_l \setdef l \in \nat\rangle$ of pairwise disjoint infinite sets such that $B'_l \subseteq B_l$.
To see this, define by recursion $B^*_0 := B_0$,
$B^*_{l+1} =  B_{l+1} \cap B^*_l$ if this set is infinite, 
and $B^*_{l+1} = B_{l+1} \setminus  B^*_l$ otherwise.
Let $D$ be the set of $l$ where the first of these two cases obtains.
As $\langle B'_l \setdef l \in D\rangle$ is $\subseteq$-decreasing, we can find $B^*_\infty$ such that for each $l \in D$, $B^*_\infty \subsetfin[] B^*_l$.
Fix an arbitrary partition $\langle P_l \setdef l \in D\rangle$ of $B^*_\infty$.
Let
$B'_l = P_l\cap B^*_l$ for $l \in D$ and $B'_l = B^*_l$ for $l \in \nat \setminus D$.
By construction, $\langle B'_l \setdef l \in \nat\rangle$ is a sequence of disjoint infinite sets and $B'_l \subseteq B_l$.
Finally, letting $\tD Z^l = \tD{A}^l  \cap (B'_l \times \nat)$, \ref{Z.refine} and \ref{Z.disjoint} obtain
by construction.
\end{proof}

Given $A \in [\nat]^{\infty}$, let us write $l \consec[A] m$ to mean that $l$ and $m$ are consecutive elements of $A$, i.e.,
\[
l \in A \text{ and } m = \min\big(A\setminus (l+1)\big).
\]

\begin{rem}\label{r.sequence}
Recall that infinite Suslin sets are Dedekind infinite in the absence even of $\CCR$ (countable choice for sets of reals) by taking left-most branches. 
Thus, a sequence $\vec z$ as above exists in \ZF{} for Suslin $\mathcal A$.
In the case of analytic $\mathcal A$ this is amounts to using 
Jankov-von Neumann, or $\sigma(\Sigma^1_1)$ uniformization for analytic relations (see~\cite[18.1]{kechris}).
If $\mathcal A$ is not assumed to be Suslin we use $\DCR$ to find $\vec z$.
 In fact, $\CCR$ suffices, as it implies that infinite sets are Dedekind infinite. 
We make another implicit appeal to $\DCR$ for the assertion that $\RUnif(\Gamma) +$ ``every set in $\Gamma$ has the Ramsey property'' implies $\RUnif^*(\Gamma)$.\footnote{We want to thank participants of the Delta Workshop in Wuhan, 2023, for the discussion regarding these points.}
For more information on $\RUnif(\Gamma)$, see \cite{muller}.
\end{rem}

\begin{dfn}[\emph{The two-dimensional tilde operator}]\label{d.2dtilde}
For each $A\in [\N]^\infty$, we define $\from[A]\in\finsqpp$ by
\[
\from[A] = \{ \hat{\tD{Z}}^l(m,n) \setdef l \consec m < n \in A\}.
\]
\end{dfn}

\medskip

{\it Remark:} 1) It is easy to see that $\from[A]\in\finsqpp$: Since $A$ is infinite, whenever $l\consec m$, there are infinitely many elements of the $m$'th vertical of $\tD{Z}^l$ that are elements of $\from[A]$; so all non-empty verticals of $\from[A]$ are infinite. Since there are infinitely many $l,m\in A$ such that $l\consec m$, we also have that $\dom(\from[A])$ is infinite. So $\from[A]\in\finsqpp$.

2) Clearly, the definition of the tilde operator depends on the sequence $\vec{z}=\langle \tD{Z}^l\setdef l\in\nat\rangle$. There will rarely be a reason to make this dependence explicit, but if there is, we'll write $\fromx[\vec{z}]{A}$.

\medskip

Having defined the tilde operator, we can now state the main technical theorem of this section:

\begin{thm}\label{t.2D.analytic}
Let $\mathcal A$ be an analytic $\finsq$-almost disjoint family. Suppose $\vec{z} = \langle \tD{Z}^l\setdef l\in\nat\rangle$
is a sequence of pairwise disjoint sets in $\finsqpp$ satisfying \ref{Z.refine} and \ref{Z.disjoint} from above.
Then the set
$$
\{W \in [\nat]^{\infty}\mid (\forall \tD{A}\in\mathcal A)\ \tD{A}\cap\fromx[\vec{z}]{W}\in \fin^2\}
$$
is Ramsey co-null. In particular, for any $W \in [\nat]^{\infty}$ there is $W_0 \in [W]^{\infty}$ such that  $\fromx[\vec{z}]{W_0}$ is $\finsq$-a.d.\  from every element of $\mathcal A$, whence $\mathcal A$ is not maximal. 
\end{thm}
 
\subsection{Properties of the tilde operator.} For use when we prove Theorem~\ref{t.2D.analytic}, we now establish some general properties of the tilde operator introduced above:
\begin{facts}~\label{f.from} 
\begin{enumerate}
\item\label{i.from.invariant} \emph{Invariance Principle:} Suppose $A,A' \subseteq \nat$.
If $A \mathbin{E_0} A'$, i.e., if 
$A \Delta A'\in \fin$, then also $\from[A] \Delta \from[A'] \in \fin[2]$.

\item\label{i.from.pigeonhole} \emph{The Pigeonhole Principle, first level:} 
For any $X \subseteq \nat$,  $A \in [\nat]^{\infty}$ and any $k\in\nat$ there is $B \in [A\cap k, A\setminus k]^{\infty}$ such that 
$\dom(\from[B])\subseteq_{\fin}  X$ or $\dom(\from[B])\subseteq_{\fin} \nat \setminus X$.

\item\label{i.from.pigeonhole.vertical} \emph{The Pigeonhole Principle, second level:} 
Suppose as above $X \subseteq \nat$,  $A \in [\nat]^{\infty}$, and $k\in\nat$ and 
$m\in\dom(\from[A])$ are given.
Then there is $B \in [A\cap k, A\setminus k]^{\infty}$ 
such that $m\in\dom(\from[B])$ and $\from[B](m) \subseteq_{\fin}X$ 
or $\from[B](m)\subseteq_{\fin} \nat\setminus X$.

\item\label{i.from.avoid.single} \emph{The Almost Disjointness Principle:} 
For $\tD A \in \mathcal A$ and any $A \in [\nat]^{\infty}$ there is $B \in [A]^{\infty}$ such that
$\from[B] \cap \tD A \in \finsq$.
\end{enumerate}
\end{facts}
\begin{proof}
(\ref{i.from.invariant}) This is because only finitely many verticals are affected when replacing $A$ with $A'$. 

\medskip
(\ref{i.from.pigeonhole})
Define a coloring
\[
c\colon [A\setminus k]^2 \to 2
\]
as follows: For $\{l,m\} \in [A\setminus k]^2$ such that $l<m$ let
\[
c(l,m)=\begin{cases}
0 &\text{ if $\widehat{\dom(\tD{Z}^l)}(m) \in X$, }\\
1 &\text{ if $\widehat{\dom(\tD{Z}^l)}(m) \notin X$.}
\end{cases}
\]
By the Infinite Ramsey's Theorem find $H\in[A\setminus k]^{\infty}$ such that $c$ takes only one color on $[H]^2$.
Letting $B=(A\cap k)\cup H$, clearly $B$ is as desired.

\medskip
(\ref{i.from.pigeonhole.vertical})
Fix $l_0,m_0$ such that 
$l_0 \consec[A] m_0$ and
$\widehat{\dom(\tD{Z}^{l_0})}(m_0)=m$.
We may assume that $k > m_0$ (otherwise increase $k$), that is,
$m \in \dom(\from[A\cap k])$.
One of the two sets
\begin{align*}
B_0&=\{n\in A\setminus k\setdef \hat{\tD{Z}}^{l_0}(m_0,n)\in X\},\\
B_1&=\{n\in A\setminus k\setdef \hat{\tD{Z}}^{l_0}(m_0,n)\notin X\}
\end{align*}
has to be infinite. Let $B=(A\cap k)\cup B_i$ where $i$ is chosen so that $B_i$ is infinite; $B$ is the desired set. 

\medskip
(\ref{i.from.avoid.single}) If for some $l\in\nat$ we have $\tD{Z}^l \subseteq \tD A$, it suffices to take $B= A$, 
since $\from[A]$ meets $\tD Z^l$ in at most one vertical.  
Otherwise, as $\tD A \in \mathcal A$, for each $l\in\nat$, $\tD A \cap \tD{Z}^l \in \finsq$.
In particular we may choose $m_l$ such that for every $m' \geq \widehat{\dom(\tD{Z}^l)}(m_l)$,
$\tD{Z}^l(m') \cap \tD A$ is finite.

It is straightforward to find a strictly increasing sequence %
$\bar m_0, \bar m_1, \hdots\in A$ such that
$(\forall k \in \nat)\:
\bar m_{k+1} >m_{\bar m_{k}}$. %
Finally let  
$B=\{\bar m_{k}\setdef k\in\nat\}$.
Then $(\from[B] \cap \tD A \cap \tD{Z}^l)(m)$ is finite for each $l,m \in B$,  so $\from[B]\cap \tD A \in \finsq$.
\end{proof}
\medskip

\subsection{A family of invariant trees in dimensions two} Let $\mathcal A$ be an analytic almost disjoint family, and fix a tree $T$ on $2\times\omega$ such that $\mathcal A=\pi[T]$. In this section, we define a map $\tD{X}\to T^{\tD{X}}$
which associates to each ${\tD{X}}\subseteq S=\N^2$ a tree $T^{\tD{X}}$ which we may think of as a ``search tree'' of attempts at finding $\tD{A}\in\mathcal A$ such that $\tD{X}\cap\tD{A}\notin\fin^2$. 
We also define a ``conditional'' version of this map, $\tD{X}\to T^{\tD{X},d}$, where the search is restricted to those attempts whose witness has infinite verticals above each $m \in d$.
A central feature of the map $\tD{X}\to T^{\tD{X}}$ is that it has good invariance properties.

\begin{dfn}
Given ${\tD{X}} \subseteq S$, define $T^{\tD{X}}$ as follows:
\[
T^{\tD{X}} = \{t \in T \setdef (\exists \tD A \in \pi[T_t]) \; \tD A \cap {\tD{X}} \notin \finsq\}.
\]
\end{dfn}

It is easy to see that $T^{\tD{X}}$ has the following properties:
\begin{facts}~\label{f.T^z}
\begin{enumerate}
\item\label{i.T^z.pruned} $T^{\tD{X}}$ is pruned, i.e., for any $t \in T^{\tD{X}}$ it holds that $[T^{\tD{X}}_t]\neq \emptyset$.
\item\label{i.T^z.non-empty} 
$(\exists \tD A \in \mathcal A)\;\tD A \cap {\tD{X}}\notin \finsq$ if and only if $\emptyset \in T^{\tD{X}}$, 
and by the previous item, this means if and only if $[T^{\tD{X}}]\neq \emptyset$.
\item\label{i.T^z.invariant} \emph{Invariance:} ${\tD{X}} \Delta {\tD{X}'}\in \finsq \Rightarrow T^{\tD{X}} = T^{{\tD{X}'}}$.
\end{enumerate}
\end{facts}

To deal with $\finsq$, we refine this definition by introducing a family of trees which also take into account some information about which verticals of the intersection are infinite:
For $\tD X \subseteq \nat^2$ and $d \in [\nat]^{<\infty}$ define
\begin{equation*}
T^{\tD X,d}=\{s \in T\setdef (\exists \tD A \in \pi[T_s])\; \tD A \cap \tD X \in\finsqp \land
 (\forall n\in d)\; \tD A(n) \cap \tD X(n) \in \finp\} 
\end{equation*}
It is again easy to see that for any $\tD X,\tD Y \subseteq \nat^2$ and $d \in [\nat]^{<\infty}$ the following hold:
\begin{facts}~\label{f.T^z,d}
\begin{enumerate}
\item $T^{\tD X,\emptyset} = T^{\tD X}$,

\item $T^{\tD X,d}$ is pruned, i.e., $t \in T^{\tD X,d} \iff [T^{\tD X,d}]\neq \emptyset$,

\item\label{i.T^z,d.subtrees} $d \subseteq d' \Rightarrow$ $T^{\tD X,d'}$ is a subtree of $T^{\tD X, d}$,

\item\label{i.T^z,d.max}
$(\exists \tD A \in \mathcal A)\;\tD A \cap \tD{X}\notin \finsq$ if and only if $[T^{\tD{X},\emptyset}]\neq \emptyset$.

\item\label{i.T^z,d.invariant} \emph{The Principle of Conditional Invariance:} 
\[
\left[(\forall n\in d)\; \tD X(n) \Delta \tD Y(n) \in \fin \land \tD X \Delta \tD Y \in\finsq \right]\Rightarrow T^{\tD X,d} = T^{\tD Y,d}.
\]

\end{enumerate}
\end{facts}
\subsection{The branch lemma in dimension two} The following lemma captures the combinatorial property of $T^{\tD{X}}$ and $\langle T^{\tD{X},d} \setdef d \in [\nat]^{<\infty} \rangle$ that is most crucial for the purpose of proving Theorem~\ref{t.2D.analytic}.
(In the way of intuition, we offer the following: The $\finsq$-almost disjointness of $\pi[T^{\tD X}] \subseteq \mathcal A$ can be isolated to a single dimension. That is, if you try to build an intersection in $\finsqp$ of two branches and $\tD X$, 
you will reach a point where  you cannot grow your intersection anymore in one of two ways: 
\eqref{i.2D.branch.I}, which expresses that you cannot grow horizontally, or \eqref{i.2D.branch.II} which expresses that you cannot grow vertically.)
We will later use this branch lemma in tandem with a ``weak genericity'' property of $\tD X$.
\begin{lem}[Branch Lemma]\label{l.2D.branch}
Suppose $\tD X \subseteq \nat^2$ and $t_0, t_1 \in T^{\tD X}$ with $\lh(t_0) = \lh(t_1)$ but $\pi(t_0)\neq \pi(t_1)$.
There is $d \in [\nat]^{<\infty}$, 
$t'_i \in  T^{\tD X,d}_{t_i}$ for each $i\in \{0,1\}$, and $m \in \nat$ such that one of the following holds:
\begin{equation*}%
\bigcap_{i \in \{0,1\}} \left\{m' \in \nat \setdef   t'_i \in T^{\tD X, d \cup \{ m'\}} \right\} \subseteq m
\tag{I$_2$}\label{i.2D.branch.I}
\end{equation*}
\begin{equation*}
(\exists m_0 \in d)
\left(\forall (t''_0, t''_1) \in  T^{\tD X, d}_{t'_0} \times T^{\tD X, d}_{t'_1}\right)\;
\big(\pi(t''_0)\cap \pi(t''_1) \cap \tD X \big)(m_0) \subseteq m 
\tag{II$_2$}\label{i.2D.branch.II}
\end{equation*}
\end{lem}
\begin{rem}
When we discuss the higher dimensional generalization in the next section, 
we will find a way to view Items \eqref{i.2D.branch.I} 
and \eqref{i.2D.branch.II} as special cases (one for each dimension) of a single statement, namely \eqref{e.T.extend}, cf.\ Fact~\ref{f.T.extend}. 
In the present case of two dimensions we do not believe that phrasing the lemma in this way would make it more transparent,
and the above formulation is the one we found most convenient.
Likewise the Pigeonhole Principles, 
Items \eqref{i.from.pigeonhole} and \eqref{i.from.pigeonhole.vertical} of Facts~\ref{f.from} above, 
become one in \ref{f.from.alphaD}\ref{i.alphaD.pidgeon}.
\end{rem}

Note that in Lemma~\ref{l.2D.branch} above, $\tD X$ is merely a parameter---it may help the reader to first set $\tD X=\nat^2$ mentally and then convince 
themselves that the proof goes through for arbitrary $\tD X$ as in the lemma (the case $\tD X\in \finsq$ is entirely uninteresting, but even then the lemma holds vacuously).
The proof idea is that the failure of the lemma would give a recipe for building two branches through $T$ whose projections have a $\finsq$-large intersection.

\medskip

Turning now to the matter of proving Lemma~\ref{l.2D.branch}, it is convenient to first define the following notation.
\begin{dfn}
Given $\tD b, \tD b' \in [\nat^2]^{<\infty}$, write $\tD b \sqsubset_2 \tD b'$ to mean that
\[
\dom(\tD b)\sqsubsetneq \dom(\tD b')
\]
 and for each $n\in \dom(\tD b)$,
$\tD b(n)\sqsubsetneq \tD b'(n)$ (where of course for $a,a'\in[\nat]^{<\infty}$,  $a \sqsubsetneq a'$
 means that
$a$ is proper initial segment of $a'$).
 \end{dfn}
 We point out the following fact, which the reader can easily prove:
 \begin{fact}
 If $\langle \tD b_k\setdef k\in \nat\rangle$ is a sequence such that $\tD b_k \sqsubset_2 \tD b_{k+1}$ for each $k\in\nat$, 
 then $\bigcup_{k\in\nat} \tD b_k \in \finsqpp$.
 \end{fact}
\begin{proof}[Proof of Lemma~\ref{l.2D.branch}]
We naturally assume $\tD X \in \finsqp$.
Fix $t_0, t_1 \in T^{\tD X}$ such that $\lh(t_0) = \lh(t_1)$ but $\pi(t_0)\neq \pi(t_1)$.
Towards a contradiction, suppose the lemma fails. 
We build sequences
\begin{itemize}
\item $\tD{b}^0, \tD{b}^1, \tD{b}^2, \hdots$ from $[\nat^2]^{<\infty}$,
\item $t^0_i, t^1_i, t^2_i \hdots$ from $T_{t_i}$ for each $i\in\{0,1\}$ 
\end{itemize}
such that %
for each $k\in \nat$
\begin{enumerate}[label=(\Alph*)]
\item\label{IH.grow}\label{IH.first}  $\tD{b}^k \sqsubset_2 \tD{b}^{k+1}$,
\item\label{IH.subseteq} $\tD{b}^{k+1} \subseteq \pi(t^k_0)\cap \pi(t^k_1)$,
\item\label{IH.keep.going} $t^{k+1}_i \in T^{\tD X,\dom(\tD{b}^k)}_{t^k_i}$.
\end{enumerate}
For $k=0$ let $\tD{b}^0=\emptyset$ and $t^0_i = t_i$ for each $i\in\{0,1\}$.
Now suppose we have constructed $\tD{b}^k$, $t^k_0$, and $t^k_1$. %

Since the lemma fails for $d=\dom(\tD{b}^k)$, $t'_i=t^k_i$ for each $i\in\{0,1\}$, and $m =  \max(d)$,
by the negation of  (\ref{i.2D.branch.I}) we can find 
$m'> m$
so that letting $d'= \dom(\tD{b}^k)\cup\{m'\}$, 
we have
\[
(t'_0, t'_1) \in  T^{\tD X,d'} \times T^{\tD X,d'}.
\] 
Now use the failure of (\ref{i.2D.branch.II}) with $d'$ substituted for $d$, finitely many (namely, $|d'|$) 
times to extend $t'_i$ to $t^{k+1}_i \in T^{\tD X,d'}$ for each $i\in\{0,1\}$ 
so that there is $\tD{b}^{k+1} \subseteq \pi(t^{k+1}_0)\cap \pi(t^{k+1}_0)$ with $m'\in \dom(\tD{b}^{k+1})$ and $\tD{b}^k\sqsubset_2 \tD{b}^{k+1}$.

Thus by our assumption that the lemma fails, we can build infinite sequences $\tD{b}^k, t^k_0, t^k_1$ for $k\in \nat$ as above;
but this contradicts that $\mathcal A$ is $\finsq$-almost disjoint since $ \bigcup_{k\in\nat}\tD{b}^k \in \finsqp$ 
and $\langle t^k_i \setdef k \in \nat\rangle$ gives us a branch $w_i\in [T_{t_i}]$ for each $i\in\{0,1\}$ such that
$\pi(w_0) \neq \pi(w_1)$ and 
\[
 \bigcup_{k\in\nat}\tD{b}^k\subseteq \pi\left(w_0\right) \cap \pi\left(w_1\right) \qedhere
\]
\end{proof}

\medskip

\subsection{The proof of Theorem~\ref{t.2D.analytic}} Having established Lemma~\ref{l.2D.branch}, we can now prove Theorem~\ref{t.2D.analytic}.
\begin{proof}[Proof of Theorem~\ref{t.2D.analytic}]
Let $\mathcal A$ and $\vec{z}$ as in the lemma be given.
Fix a tree $T$ on $2\times\nat$ such that $\pi[T]=\mathcal A$. 

The first step in the proof is Claim \ref{c.2D.constant} below, which shows that, for each $d \in [\nat]^{<\infty}$, we can make the function 
\[
B \mapsto T^{\from[B],d}
\] 
constant for a ``large'' set of $B\in [\nat]^\infty$, by exploiting the invariance properties of the trees we established above, and the fact that analytic sets are completely Ramsey.
The following claim can also be interpreted as stating the existence of a ``weakly generic'' set $W_0$.
It is crucial here that the branch lemma is a statement about countably many countable trees.

\begin{claim}\label{c.2D.constant}
For every $W \in [\nat]^{\infty}$ there is $W_0 \in [W]^{\infty}$ such that 
\[
(\forall W' \in [W_0]^{\infty})(\forall d \in [\dom(\from[W'])]^{<\infty})\; T^{\from[W'],d}=T^{\from[W_0],d}
\]
\end{claim}
\begin{proof}[Proof of Claim~\ref{c.2D.constant}]
The proof is by a fusion argument.
Let 
\[
\langle(d_k, t_k)\setdef k\in\nat\rangle
\]
 enumerate $[\nat]^{<\infty}\times T$ so that each element of $[\nat]^{<\infty}\times T$ occurs infinitely many times in the enumeration.
Let for each $k\in \nat$
\[
D_k=\{B \in [\nat]^\infty\setdef t_k \in T^{\from[B],d_k}\}
\]
and note that this set is analytic and hence, completely Ramsey.

We inductively build a sequence $(b_k,C_k) \in [\nat]^{<\infty}\times[\nat]^{\infty}$ for $k\in \nat$.
Let $(b_0, C_0)=(\emptyset, W)$. Suppose now we have constructed $(b_k,C_k)$.
Find $C_{k+1} \in [C_k\setminus b_k]^{\infty}$ such that
\begin{equation}\label{e.make.constant.one.step}
[b_k, C_{k+1}]^{\infty}\subseteq D_k \text{ or } [b_k, C_{k+1}]^{\infty}\subseteq [\nat]^{\infty} \setminus D_k,
\end{equation}
and let $b_{k+1}=b_k\cup\{\min C_{k+1}\}$.
Finally let $W_0 = \bigcup_{k\in\nat} b_k$. 

To see that $W_0$ is as desired, let an arbitrary $W'\in[W_0]^{\infty}$ and a finite set $d \subseteq \dom(\from[W'])$ be given.
To show $T^{\from[W'],d}=T^{\from[W_0],d}$, consider an arbitrary $t \in T$.
Find $k\in\nat$ such that $(d_k,t_k)=(d,t)$ and $d \subseteq \dom(\from[b_k])$.
Since $W' \subseteq_{\fin} C_{k+1}$ we may choose $W''$ such that $W'' \mathbin{E_0} W'$ and $W'' \in [b_k, C_{k+1}]^{\infty}$.
Thus by Equation~\eqref{e.make.constant.one.step},
\[
t \in T^{\from[W''],d} \iff t \in T^{\from[C],d}, 
\]
with $C = d_k \cup C_{k+1}$.
By the conditional invariance of this tree expressed in Fact\ref{f.T^z,d}\eqref{i.T^z,d.invariant} and since
$\from[W''] \Delta \from[W'] \in \finsq$ and $d \subseteq \dom_\infty(\from[W''])$ and so also $\from[W''](n) \Delta \from[W'](n) \in \fin$ for each $n\in d$, 
we also have 
\begin{equation}\label{e.fusion}
t \in T^{\from[W'],d} \iff t \in T^{\from[C],d}.
\end{equation}
Since in particular, \eqref{e.fusion} also holds with $W'$ replaced $W_0$, and since $t\in T$ was arbitrary, we conclude $T^{\from[W_0],d}=T^{\from[W'],d}$.
\renewcommand{\qedsymbol}{{\tiny  Claim~\ref{c.2D.constant}~}$\Box$}
\end{proof}

Now fix $W_0$ as in the previous claim, and let $T^{*,d}=T^{\from[W_0],d}$; then $T^{\from[W],d}=T^{*,d}$ for all $W\in [W_0]^\infty$ with $d \in [\dom(\from[W])]^{<\infty}$. 
Intuitively, $W_0$ is ``generic enough'' to ensure that $\from[W_0]$ will be disjoint from every set from $\mathcal A$. The following two claims each use this ``weak genericity''; the first does so in a particularly simple manner.

\begin{claim}\label{c.2D.not1branch}
$\lvert \pi[T^{*,\emptyset}]\rvert \neq 1$.
\end{claim}
\begin{proof}[Proof of Claim~\ref{c.2D.not1branch}]
Towards a contradiction suppose 
\[
\{\tD A\} = \pi[T^{*,\emptyset}].
\]
Find $W_1 \in [W_0]^{\infty}$ such that $\from[W_1] \cap \tD A \in \finsq$, using Fact~\ref{f.from}(\ref{i.from.avoid.single}).
Since $T^{\from[W_1]} = T^{*,\emptyset}$ by the previous claim, we have  $\{\tD A\} = \pi[T^{\from[W_1]}]$. This leads to a contradiction as follows: Since $\pi[T^{\from[W_1]}]\neq \emptyset$, we have in particular $\emptyset \in T^{\from[W_1]}$, so there is $w \in [T]$ with $\pi(w) \cap \from[W_1] \in \finsqp$.
But clearly also $w \in [T^{\from[W_1]}]$ and so $\pi(w) = \tD A$. Thus $\tD A \cap \from[W_1] \in \finsqp$, contradicting that we chose $W_1$ such that $\tD A \cap \from[W_1] \in \finsq$.
\renewcommand{\qedsymbol}{{\tiny  Claim~\ref{c.2D.not1branch}~}$\Box$}
\end{proof}
We continue exploring the ``weak genericity'' of $W_0$, this time using the pigeonhole principles proved in Fact~\ref{f.from} as well as Lemma~\ref{l.2D.branch}.
\begin{claim}\label{c.2D.one.branch}
$\pi[T^{*,\emptyset}]=\emptyset$.
\end{claim}
\begin{proof}[Proof of Claim~\ref{c.2D.one.branch}]
Otherwise by the previous claim, $\lvert \pi[T^{*,\emptyset}]\rvert\geq 2$.
Then we can choose $t_0, t_1 \in T^{*,\emptyset}$ with $\lh(t_0) = \lh(t_1)$ and $\pi(t_0)\neq \pi(t_1)$.
Use Lemma~\ref{l.2D.branch} to find $d \in [\nat]^{<\infty}$ and 
$t'_i \in  T^{*,d}_{t_i}$ for each $i\in \{0,1\}$ satisfying Clause
(\ref{i.2D.branch.I}) or (\ref{i.2D.branch.II}) of said lemma.
We now argue by cases.

\medskip
\paragraph*{First case: Clause (\ref{i.2D.branch.I}) of Lemma~\ref{l.2D.branch} holds.} 
For each $i\in\{0,1\}$ let
\[
X_i = \bigcup\{\dom_\infty\big(\pi(w)\cap\from[W_0]\big) \setdef w \in [T^{*,d}_{t'_i}] .
\]
Then $X_0 \cap X_1 \in \fin$:
Otherwise, there are unboundedly many $m' \in X_0 \cap X_1$. 
For each such $m'$,  
there is $(w_0, w_1) \in T^{*,d}_{t'_0} \times T^{*,d}_{t'_1}$
such that 
\[
m \in \bigcap_{i\in \{0,1\}} \dom_\infty(\pi(w_i) \cap \from[W_0].
\]
But then by definition, $w_i$ witnesses that $t'_i  \in  T^{*,d\cup \{m'\}}_{t'_i}$, for each $i \in \{0,1\}$.
Since this happens for unboundedly many $m'$, we contradict (\ref{i.2D.branch.I}).

We have shown $X_0 \cap X_1 \in \fin$.
Use Facts~\ref{f.from}(\ref{i.from.pigeonhole}) to find $W_1 \in [W_0]^{\infty}$ such that $d \subseteq \dom(\from[W_1])$ and such that
$\dom(\from[W_1])\subseteq_{\fin} X_0$ or $\dom(\from[W_1])\subseteq_{\fin} \nat \setminus  X_0$.
Thus, we may fix $i\in\{0,1\}$ so that 
\begin{equation}\label{e.2D.small}
\dom(\from[W_1])\cap X_i \in \fin.
\end{equation}

But by the definition of the invariant tree, this contradicts that  $t'_i \in T^{*,d}$ as we now show in detail for the reader's convenience:
Since $T^{*,d} = T^{\from[W_1],d}$, $t'_i  \in T^{\from[W_1],d}$ and so by definition of the tree there is $w \in [T_{t'_i}]$ such that 
\begin{equation}\label{e.2D.big}
\pi(w) \cap \from[W_1] \in \finsqp
\end{equation} 
and
$d \subseteq \dom_\infty(\pi(w) \cap \from[W_1])$. 
Also by definition of the tree, 
$w$  is a branch through $T^{W_1,d}_{t'_i}$ and hence through $T^{*,d}_{t'_i}$. 
But then $\dom_\infty\big(\pi(w) \cap \from[W_1]\big) \subseteq X_i$ 
which using \eqref{e.2D.big} entails $\dom(\from[W_1]) \cap X_i \in \finp$,
contradicting Equation \eqref{e.2D.small}.

\medskip
\paragraph*{Second case: Clause (\ref{i.2D.branch.II}) of Lemma~\ref{l.2D.branch} holds.}
Fix $m\in d$ as in Clause (\ref{i.2D.branch.II}). For each $i\in\{0,1\}$ let
\[
X_i = \bigcup\big\{\pi(w)(m) \setdef w \in [T^{*,d}_{t'_i}]\big\}.
\]
By assumption $X_0 \cap X_1 \in \fin$.
Argue as in the previous case, this time using Facts~\ref{f.from}(\ref{i.from.pigeonhole.vertical}) to find $W_1 \in [W_0]^{\infty}$ such that $d \subseteq \dom(\from[W_1])$ and such that we may find $i\in\{0,1\}$ with
\begin{equation}\label{1D.case.B.W_1.smallintersection}
\from[W_1](m)\cap X_i(m) \in \fin.
\end{equation} 
Then similarly to the previous case, this contradicts that $m\in d$ and
$t'_i \in T^{\tilde W_1,d}$.
Having reached a contradiction in each case, this proves the claim.
\renewcommand{\qedsymbol}{{\tiny  Claim~\ref{c.2D.one.branch}~}$\Box$}
\end{proof}

Since $T^{*,\emptyset}=T^{\from[W_0],\emptyset}$, the last two claims show that $[T^{\from[W_0],\emptyset}]=\emptyset$. 
By Facts~\ref{f.T^z,d}(\ref{i.T^z,d.max}) we conclude that $\mathcal A$ is not maximal.
\renewcommand{\qedsymbol}{{\tiny  Theorem~\ref{t.2D.analytic}~}$\Box$}
\end{proof}

\subsection{Proof of the main theorem in two dimensions}\label{proof.t.2D.precise}

Now we can easily prove Theorem~\ref{t.2D.precise}. 
\begin{proof}[Proof of theorem~\ref{t.2D.precise}]
Suppose to the contrary $\mathcal A$ is an infinite $\finsq$-\emph{mad} family in $\Gamma$.
Using $\DCR$ (or its consequence, $\CCR$) choose a sequence $\langle \tD{A}^l\mid l\in\N\rangle$ of elements in $\mathcal A$, and let $\vec{z}=\langle \tD{Z}^l\setdef l\in\nat\rangle$ 
be defined from $\langle \tD{A}^l\mid l\in\N\rangle$ as in Lemma \ref{l.Z.2D}.
We reduce the proof to the case that $\mathcal A$ is analytic:
\begin{claim}\label{c.make.analytic}
There is $W_1\in [\nat]^{\infty}$ and an analytic set $\mathcal A'$ such that 
for any $W \in [W_1]^{\infty}$ there is $\tD A \in \mathcal A'$ with 
$\fromx[\vec{z}]{W} \cap \tD A \notin \finsq$.
\end{claim} 
\begin{proof}[Proof of Claim~\ref{c.make.analytic}]
Let $R(W,\tD A)$ be the relation defined by
\[
R(W,\tD A) \iff W \in [\nat]^{\infty} \land \tD A \in \mathcal A \land \fromx[\vec{z}]{W}\cap \tD A \notin \finsq.
\]
By our assumptions on $\Gamma$, $R$ is in $\Gamma$.
So we may find $W_0 \in [\nat]^{\infty}$ and $f$ in $\Gamma'$ such that
for all $W\in [W_0]^{\infty}$, $R(W, f(W))$.

Recalling from \S\ref{s.prelim} that we identify $\powerset(\nat^2)$ with a subset of $2^\nat$ 
and denoting by $N_s$ the basic open neighborhood in 
$2^\nat$ given by $s\in 2^{<\nat}$, for each $s\in2^{<\nat}$, the set
\[
\{W\in [W_0]^{\infty} \setdef f(W) \in N_s\}
\]
is in completely Ramsey. 
An easy diagonalization shows that we may thin out $W_0$ to $W_1 \in [W_0]^{\infty}$ 
so that $f$ restricted to $[W_1]^{\infty}$ is continuous. 
Letting $\mathcal A' = \ran(f)$ this proves the claim.
\renewcommand{\qedsymbol}{{\tiny  Claim~\ref{c.make.analytic}~}$\Box$}
\end{proof}
\noindent By Theorem~\ref{t.2D.analytic} no such $\mathcal A'$ can exist, finishing the proof.
\renewcommand{\qedsymbol}{{\tiny  Theorem~\ref{t.2D.precise}~}$\Box$}
\end{proof}

\section{The general case}\label{s.alphaD}

We now treat the general case of ideals arising from $\fin$ by iterated Fubini products, 
that is, ideals from $\mathfrak F$. 
This class contains the class of iterated Fréchet ideals as defined in \cite{jackson-kechris-louveau}.
We show:

\begin{thm}[\DCR]\label{t.Finalph}
Suppose $\Gamma$ is a pointclass such that $\RUnif^*(\Gamma)$ holds and $\mathcal I \in \mathfrak F$.
Then there are no infinite $\mathcal I$-\emph{mad} families in $\Gamma$.
\end{thm}

\subsection{The ideal}
We need a convenient representation of an arbitrary ideal $\mathcal I \in \mathfrak F$.
We start with the domain:
The easiest case is $\fin[k]$ for $k < \omega$,
whose domain will be $\nat^k$, which we identify with the set of sequences from $\nat$ of length $k$.
More generally, we will think of the domains of our ideals in terms of leaves of trees.
This needs some notational
preparations, which we now undertake.

Let $\tD X \subseteq \nat^{<\nat} $ consist of pairwise $\subseteq$-incompatible sequences, that is,
if $s, s' \in \tD X$ and $s \neq s'$, then the longest common intial segment of $s$ and $s'$ is $\emptyset$, the empty sequence and root of the tree $\nat^{<\nat}$.
Define\footnote{This notation is of course inspired by the two-dimensional case, 
where we chose to view subsets of $S = \nat^2$ as binary relations on $\nat$ instead of trees.} 
\[
\higherdom(\tD X)=\{s\res k \setdef s \in \tD X, k < \lh(s)\}.
\]
Note that $\higherdom(\tD X) \cup \tD X$ is a subtree of $\nat^{<\nat} $, ordered by $\subseteq$,
and $\tD X$ consists of this tree's maximal elements, i.e., its leaves.
Given $s, s' \in \nat^{<\nat}$, 
we write $s \conc s'$ for the sequence $t$
such that $t \res \lh(s) = s$ and
$t(k) = s'\big(k - \lh(s)\big)$ for $k < \lh(s) + \lh(s')$,
and $s \conc n$ for $s \conc \langle n \rangle$ when $n \in \nat$. %
For $s \in \nat^{<\nat}$
define
\[
\splitset(\tD X,s) = \{n\in\nat\setdef(\exists s' \in\nat^{<\nat} )\; s \conc  n  \conc s' \in \tD X\},
\]
the splitting set of $\tD X$ at $s$.
Provided that $n< \lvert \splitset(\tD X,s) \rvert$, we shall also write
\[
\splitenum(\tD X,s,n) := \widehat{\splitset(\tD X, s)}(n)
\]
as a short-hand for the $n$th element in the increasing enumeration of  $\splitset(\tD X,s)$.
Moreover, define 
\[
\tD X(s) = \{s'\in\nat^{<\nat} \setdef s \conc s' \in \tD X\}.
\]

We are now ready to describe our presentation of the ideals in $\mathfrak F$.
Let $U$ be a well-founded subtree of $\nat^{<\nat}$,
and assume that for each $s \in U$, 
$\splitset(U,s) = \nat$ or $s$ is a leaf of $U$. 

For such $U$, we let $S_U$ be the set of leaves of $U$.
We define the ideal $\fin[U]$ on $S_U$ by induction on tree rank.
Two different trees might well define the same ideal.

For $U = \nat^{\leq 1}$, representing the least-ranked tree that we consider, take $\fin[U]$ to be $\fin$, which is clearly (up to some trivial rewriting) an ideal on $S_U = \nat^1$, which indeed is the set of leaves of $U$.
For $U$ of greater rank, write $S_U$ for its set of leaves and define
$\tD X \subseteq S_U$ to be an element of $\fin[U]$ if and only if 
\[
\{n\in\nat\setdef \tD X(s) \notin \fin[U(s)], \text{ where $s = \langle n\rangle$}\} \in \fin.
\]
It is clear from the definition that 
firstly,
\[
\fin[U] = \fin \bigotimes \left(\fin[U(\langle n \rangle)]\right)_{n\in\nat}
\]
and that secondly, every ideal from $\mathfrak F$ is of the form $\fin[U]$ for some $U$ as above.

\medskip

We now introduce an alternative notation. This is the notation we will use for the remainder of this article.
Instead of parametrizing our ideals by trees $U$, 
we will parametrize them by rank functions.
By a rank function we mean a partial function 
\[
\laddermap\colon\nat^{<\nat} \rightharpoonup \aleph_1
\]
such 
that for some tree $U$ as above, $\laddermap$ is the rank function on $U$, that is,
$\dom(\laddermap) = U$ and 
$\laddermap(s)$ is the rank of $U(s)$.

We write $S_{\laddermap}$ for $S_U$ and $S_{\ladderres}$ for $S_{U(s)}$.
Likewise, we write
$\fin[\laddermap]$ for $\fin[U]$ and $\fin[\ladderres]$ for $\fin[U(s)]$.
We also write $[S_{\laddermap}]^{+}_{\finalph}$ 
or $\finalphp$ for the co-ideal $\powerset(S_{\laddermap})\setminus\finalph$.

Note that $s \in S_{\laddermap}$ if and only if $\laddermap(s)= 0$, i.e., $s$ is a leaf of $U = \dom(\laddermap)$.
Moreover, for example, $\laddermap(\emptyset)$ is the rank of $U$ and $\laddermap\big(s\res (\lh(s)-1)\big) = 1$ 
for each $s \in S_{\laddermap}$.

\medskip

We will mostly be concerned with the following subset of $\finalphp$:
Define 
\[
\finalphpp = \{\tD X \in \powerset(S_\laddermap) \setdef \big(\forall s \in \higherdom(\tD X)\big) \: \splitset(\tD X,s) \in \finp\}.
\]
In other words, $\tD X \in \finalphpp$ if and only if for each $s\in \nat^{<\nat} $, $\splitset(\tD X,s)$ is infinite or empty;
the first alternative obtains if and only if $s \in \higherdom(\tD X)$.
We also denote $\finalphpp \cap \powerset(\tD X)$ by $[\tD X]^{++}_{\finalph}$.

It is straightforward to check that for $\tD X \subseteq S_\laddermap$, $\tD X \in [S_\laddermap]^{+}_{\finalph}$ 
if and only if there is $\tD X' \subseteq \tD X$ such that $\tD X' \in [S_\laddermap]^{++}_{\finalph}$. In fact, there 
is a maximal (with respect to $\subseteq$) such set $\tD X'$ which we denote by $\tD{X}^{++}$, and 
\[
\tD{X}^{++} =\{s\in \tD X\setdef \big(\forall k < \lh(s)\big)\;\tD X(s\res k)\in \finp[{\ladderres[(s\res k)]}]\}. 
\]

\medskip

\subsection{The higher tilde operator}

\begin{dfn}[The higher dimensional $\sim$-operator]
First, for readers who elected to skip the previous section, let us recall again the relevant notation introduced there:

Given $A \in [\nat]^{\infty}$, we write $l \consec[A] m$ to mean that $l$ and $m$ are consecutive elements of $A$, i.e.,
\[
l \in A \text{ and } m = \min A\setminus l+1.
\]
Recall that we use the notation $\widehat X$
 to denote the increasing enumeration $\widehat{X}\colon \lvert X \rvert \to X$
of a finite or infinite set $X \subseteq \nat$. 
We shall presently define what can be viewed as a generalization of the increasing enumeration, 
for sets
$\tD{X} \subseteq S_\laddermap$.
Somewhat vaguely, we first define $\whathat{\tD{X}}$ as the (partial) function which takes a sequence of numbers to a node in the tree $\higherdom(\tD{X})\cup \tD{X}$ using the increasing 
enumeration of each level of the tree. 
Then the partial function $\widehat{\tD X}$ is simply a restriction of $\whathat{\tD{X}}$ which is onto the set of leaves.

To be precise:
Given $\tD X \subseteq S_\laddermap$, we define a partial function
\[
\widehat{\tD{X}} \colon \nat^{<\nat}  \rightharpoonup \tD X
\]
as well as an auxiliary partial function 
\[
\whathat{\tD{X}} \colon  \nat^{<\nat}  \rightharpoonup \higherdom(\tD X) \cup \tD X 
\]
where $\whathat{\tD{X}}(s)$ is defined by induction on $\lh(s)$.
We let
\[
\whathat{\tD{X}}(\emptyset) = \emptyset.
\]
Now for $s \in  \nat^{<\nat} $ and $n\in\nat$, 
towards defining $\whathat{\tD{X}}(s\conc n)$
let us suppose by induction that $\bar s := \whathat{\tD{X}}(s)$ is already defined and
that $\bar s \in \higherdom(\tD X)\cup \tD X$.
We let
\[
\whathat{\tD{X}}(s\conc n)  = \bar s \conc
\splitenum(\tD X, \bar s,n) 
\] 
if the expression on the right is defined, with the understanding that
\[
\whathat{\tD{X}}(s\conc n)  \uparrow 
\]
if already $\bar s \in \tD X$ or if $\bar s \in \higherdom(\tD X)$ but $\lvert \splitset(\tD X, \bar s) \rvert \leq n$.
 Note that the second case never occurs if $\tD X \in \finpp[\laddermap]$.

Having thus defined $\whathat{\tD{X}}$
we can now define $\widehat{\tD{X}}$ as the restriction
\[
\widehat{\tD{X}} := \whathat{\tD{X}} \res \{s \in  \nat^{<\nat}  \setdef \whathat{\tD{X}} (s) \in \tD X\}.
\]
\begin{rem}
In the previous section, we've used the
short-hand
$\hat{\tD Z}^l(m, n)$
for the $n$th pair in the $m$th vertical of $\tD{Z}^l$.
This is consistent with our new notation:
\begin{equation*}
\hat{\tD Z}^l(m,n) = \widehat{\tD{Z}^l}(\langle m,n\rangle)\\
\end{equation*}
when $\finalph = \fin^2$, $S_\laddermap = \nat^2$.
In this sense, our notation generalizes that of the previous section.
\end{rem}

\medskip

Fix a $\finalph$-a.d.\  family $\mathcal A$.
The higher tilde operator is (again) defined relative to an
appropriate sequence $\langle \tD{Z}^l \setdef  l\in\nat\rangle$.
We require that this sequence satisfy the following:
\begin{enumerate}[label=\textup{(\Alph*$_\laddermap$)}]
\item \label{Z.refine.alph} For each $l\in \nat$, $\tD{Z}^l \in [\tD A]^{++}_{\finalph}$ for some $\tD A \in \mathcal A$.
\item \label{Z.disjoint.alph} For each $m \in \nat$ there is at most one $l \in \nat$ such that $\tD{Z}^l (\langle m \rangle) \neq\emptyset$.
\end{enumerate}
Note that as in the two-dimensional case, by \ref{Z.disjoint.alph} we also have $\tD{Z}^l \cap \tD{Z}^{l'} =\emptyset$ for $l \neq l'$.

Before defining the higher tilde operator, let us quickly verify that such a sequence is indeed available to us.
\begin{lem}\label{l.Z.alphaD}
If there is an infinite $\finalph$-\emph{mad} family, there is a sequence $\langle \tD{Z}^l \setdef  l\in\nat\rangle$ satisfying \ref{Z.refine.alph} and \ref{Z.disjoint.alph} above.
\end{lem}
The proof is identical to the proof of Lemma~\ref{l.Z.2D} in the two-dimensional case, the only differences being in notation.
\begin{proof}
Suppose $\mathcal A$ is an infinite $\finalph$-\emph{mad}.
Using $\DCR$ choose a sequence $\langle \tD{A}^l\mid l\in\N\rangle$ of elements of $\mathcal A$ (for the current lemma, $\CCR$ suffices)
and let $B_l = \splitset(\tD{A}^l, \emptyset)$.
Exactly as in the proof of Lemma~\ref{l.Z.2D}, find a sequence $\langle B'_l \setdef l \in \nat\rangle$ of pairwise disjoint infinite subsets of $\nat$ such that $B'_l \subseteq B_l$.
Letting 
\[
\tD Z^l = \{s \in (\tD A^l)^{++} \setdef s(0) \in B'_l \}, 
\]
\ref{Z.refine.alph} and \ref{Z.disjoint.alph} obtain
by construction.
\end{proof}

We now have all the prerequisites to define the higher tilde operator.
Given $A\subseteq \nat$, let $\from[A] \subseteq S_\laddermap$ be defined as follows: 
\[
\from[A] = \{\widehat{\tD{Z}^l}(\hat a) \setdef  \{l\}\cup a \in [A]^{<\infty}\wedge l \consec \min(a) \wedge \hat a \in \dom(\widehat{\tD{Z}^l})\}.
\]
\end{dfn}

We will see that $\from[A]$ behaves similarly to the two-dimensional and one-dimensional cases 
(for the latter, see \cite{pnas}). %
In particular, analogously to Facts~\ref{f.from} we have the following:
\begin{facts}\label{f.from.alphaD}~
\begin{enumerate}[label=(\arabic*),ref=(\arabic*)]
\item\label{i.alphaD.invariance} 
\emph{Invariance:}
Suppose $A,B\in [\nat]^\infty$. If $A \mathbin{E_0} B$ and 
\[
s \in \higherdom(\from[A]) \cap \higherdom(\from[B]), 
\]
then $\from[A](s) \Delta \from[B](s) \in \fin[\ladderres]$.

\item\label{i.alphaD.pidgeon} 
\emph{The Pigeonhole Principle in higher dimensions:}
For any $A\in[\nat]^{\infty}$, 
$k\in\nat$, 
$X \subseteq \nat$, and
$s \in \higherdom(\from[A])$ 
there is $A' \in [A\cap k, A\setminus k]^{\infty}$ with 
$s \in \higherdom(\from[A'])$ and such that 
\[
\splitset(\from[A'],s) \subseteq_{\fin} X \text{ or }
\splitset(\from[A'],s) \subseteq_{\fin} \nat\setminus X.
\]

\item\label{i.alphaD.tilde.ad}
\emph{The Almost Disjointness Principle:}
For any $\tD A \in \mathcal A$ and $B \in [\nat]^{\infty}$ there is
$B' \in [B]^{\infty}$ such that $\tD A \cap \from[B']\in \finalph$.
\end{enumerate}
\end{facts}
\begin{proof}
\ref{i.alphaD.invariance}
This is an easy consequence of the definition of the $\sim$-operator. 
As the reader is invited to verify,  $s \in \higherdom(\from[A])$ holds 
if and only if 
there exists $l_A \in A$ and $c_A \in  [A]^{<\infty}$ such that $l_A \consec[A] \min(c_A)$ 
and 
\[
s =  \whathat{\tD{Z}^{l_A}(\widehat{c_A})}.
\] 
Let us fix such $l_A, c_A$ and let us also fix $l_B, c_B$ witnessing $s \in \higherdom(\from[B])$ in the analogous manner. 

Find $m\in \nat$ such that  $A \setminus m = B\setminus m$ and $c_A \cup c_B \subseteq m$.
Write $C_A$ for $\{l_A\}\cup c_A \cup (A\setminus m)$
and $C_B$ for $\{l_B\}\cup c_B \cup (B\setminus m)$.
Then clearly for $n \in\nat\setminus m$ it holds that 
\[
\from[A](s\conc n)=\from[C_A](s\conc n)=\from[C_B](s\conc n) = \from[B](s\conc n).
\]
Thus 
$\splitset(\from[A] \Delta \from[B], s)$ must be finite, and hence
$\from[A](s) \Delta \from[B](s) \in \fin[\ladderres]$.

\ref{i.alphaD.pidgeon}
Let $A$, $k$, $X$, and $s$ be given. First assume $s \neq \emptyset$:
In this case, fix $l \in \nat$ and $a \subseteq A$ finite such that 
$l \consec[A] \min(a)$ and 
\[
\whathat{\tD{Z}^l}(\hat a) = s. 
\]
Increasing $k$ if necessary, we may 
assume that $k$ is large enough so that $\{l\}\cup a \subseteq k$.
Let
\begin{align*}
B_0 &= \{n \in A \setdef \splitenum(\tD{Z}^l, s, n) \in X\},\\
B_1 &= \{n \in A \setdef \splitenum(\tD{Z}^l, s, n) \notin X\}, 
\intertext{and}
A_i &= \{l\} \cup a \cup B_i,
\end{align*}
for each $i\in\{0,1\}$.
It is straightforward from the definitions that
\begin{align*}
\splitset\big(\from[A_0], s\big) &\subseteq X, \\
\splitset\big(\from[A_1], s\big) &\subseteq \nat\setminus X.
\end{align*}
Fix $i\in\{0,1\}$ so that $B_i$ is infinite and let
\[
A' = (A \cap k)\cup B_i.
\]
Since $\splitset\big(\from[A'],s\big) = \splitset\big(\from[A_i],s\big)$, $A'$
is as required.

In case $s = \emptyset$, the proof differs from the proof of Facts~\ref{f.from}, Item~\eqref{i.from.pigeonhole} only in notation:
Define a coloring
\[
c\colon [A\setminus k]^2 \to 2
\]
as follows: For $\{l,m\} \in [A\setminus k]^2$ such that $l<m$ let
\[
c(l,m)=\begin{cases}
0 &\text{ if $\splitenum(\tD{Z}^l, \emptyset, m) \in X$, }\\
1 &\text{ if $\splitenum(\tD{Z}^l, \emptyset, m) \notin X$.}
\end{cases}
\]
As previously, by the Infinite Ramsey's Theorem find $H\in[A\setminus k]^{\infty}$ such that $c$ takes only one color on $[H]^2$.
As previously, letting $B=(A\cap k)\cup H$, we obtain a set $B$ as desired.

\ref{i.alphaD.tilde.ad}
We construct a sequence 
$\langle n_k \setdef k\in\nat\rangle$ from $\nat$ by induction. 
Let $n_0 = \min(B)$.
Now suppose we have already defined $n_k$.
Since $\tD{Z}^{n_k} \cap \tD A \in \finalph$, it will be true for any large enough $m \in \nat$ that
with
 \[
 s=\whathat{\tD{Z}^{n_k}}(\langle m\rangle), 
 \]
 we have 
$(\tD{Z}^{n_k} \cap \tD A) (s) \in \fin[\ladderres]$.
Choose $n_{k+1}$ to be any number in $B$ so that $n_{k+1} > n_k$ and $n_{k+1}$ has the same property as $m$ above.
Finally, let $B' = \{n_k \setdef k\in\nat\}$.

By construction, given any two consecutive elements $l,m$ of $B'$ it holds that $(\tD{Z}^l \cap \tD A)(s) \in \fin[\ladderres]$ 
where
\[
s=\whathat{\tD{Z}^l}(\langle m\rangle).
\]
But by the definition of the higher tilde operator, every element of
$\from[B'] \cap \tD{Z}^l$ extends $s$. 
It follows that for any $s \in \nat^1$ we have
$(\from[B'] \cap \tD A)(s) \in \fin[\ladderres]$ and so also
$\from[B'] \cap \tD A \in \finalph$.
\end{proof}

\begin{rem}
In Facts~\ref{f.from.alphaD}\ref{i.alphaD.invariance} above, 
it is assumed that $s \in \higherdom(\from[A]) \cap \higherdom(\from[B])$.
It is indeed necessary to ask $s \notin \higherdom(\from[A]) \Delta \higherdom(\from[B])$; 
for otherwise, it may be the case 
that, e.g., $s \notin \higherdom(\from[A])$ and thus $\from[A](s) = \emptyset$, 
while $s \in \higherdom(\from[B])$ and thus $\from[B](s) \in \finpp[\ladderres]$; so clearly, $\from[A](s)\Delta\from[B](s)\notin \finp[\ladderres]$ in this case. 
If $s \notin \higherdom(\from[A]) \cup \higherdom(\from[B])$ it holds that $\from[A](s) = \from[B](s)=\emptyset$, 
but we do not seem to ever need this case.

In Facts~\ref{f.from.alphaD}\ref{i.alphaD.tilde.ad} above, 
we could even demand that for for any $B'' \in [B']^\infty$, $\tD A \cap \from[B''] \in \finalph$; 
we leave details to the reader, as we shall not need this stronger statement.
\end{rem}

\subsection{The higher dimensional tree family}

\begin{dfn}
For $\tD X \in \finalphp$ and $\tD d \in \Big[\higherdom(S_\laddermap)\Big]^{<\infty}$  
define
\begin{multline*}
T^{\tD X,\tD d}=\Big\{t \in T\setdef (\exists \tD A \in \pi[T_t])\; \tD A \cap \tD X \in\finalphp \;\land (\forall s\in d)\;\\
 \tD A(s) \cap \tD X(s) \in \finp[\ladderres]\Big\} 
\end{multline*}
\end{dfn}

It is again easy to see that for any $\tD X,\tD Y \in \finsqp$ and $\tD d \in \Big[\higherdom(S_\laddermap)\Big]^{<\infty}$
 the following hold:
\begin{facts}~\label{f.alphaD.T^z,d}
\begin{enumerate}

\item\label{i.alphaD.T^z,d.pruned} $T^{\tD X,\tD d}$ is pruned, i.e., $t \in T^{\tD X,\tD d} \iff [T^{\tD X,\tD d}_t]\neq \emptyset$,

\item\label{i.alphaD.T^z,d.non-empty} $T^{\tD X,\emptyset}\neq\emptyset$ 
if and only if there is $\tD A \in \mathcal A$ such that $\tD A \cap \tD X \in \finalphp$,
\item\label{i.alphaD.T^z,d.invariant} \emph{The Principle of Conditional Invariance:}
\[
[(\forall s\in \tD d)\; \tD X(s) \Delta \tD Y(s) \in \fin[\ladderres] \land \tD X \Delta \tD Y \in\finalph ]
\Rightarrow T^{\tD X,\tD d} = T^{\tD Y,\tD d}
\] 
\item\label{i.alphaD.T^z,d.invariant.corollary} 
In particular, if $W_0, W_1 \in \finp$, $W_0 \mathbin{E_0} W_1$ and for each $i\in\{0,1\}$ 
it holds that $\tD d \subseteq \higherdom(\from[W_i])$, 
then also $T^{\from[W_0],\tD d} = T^{\from[W_1],\tD d}$,
\end{enumerate}
\end{facts}
The most crucial fact about these trees is expressed by the following lemma:
\begin{lem}\label{l.alphaD.branch}\label{l.alphaD.branch'}
Suppose $t_0, t_1 \in T^{\tD X,\emptyset}$, $\lh(t_0) = \lh(t_1)$ but $\pi(t_0)\neq \pi(t_1)$. 
Then we can find
\begin{itemize} 
\item a subtree $\tD d$ of $\higherdom(S_\laddermap)$,
\item nodes  
$t'_0, t'_1 \in  T$ such that for each $i\in \{0,1\}$,
\[
t'_i \in T^{\tD X,\tD d}_{t_i}
\]
\item as well as $s \in \tD d$ and $m\in\nat$,
\end{itemize}
such that the following condition (which we denote by $\cond$, or in more detail, by $\condargs$  for later reference) is met:
if $\laddermap(s) > 1$, then
\begin{equation*}%
\label{i.main}
\bigcap_{i \in \{0,1\}}\left\{m' \in \nat \setdef   t'_i \in T^{\tD X, \tD d \cup \{ s \conc m'\}} \right\} \subseteq m
\tag{I$_\laddermap$}
\end{equation*}
and if $\laddermap(s) =1$, then
\begin{equation*}\label{i.main.cond}
\left(\forall (t''_0, t''_1) \in  T^{\tD X,\tD d}_{t'_0} \times T^{\tD X,\tD d}_{t'_1}\right)\;
\splitset\Big(\pi(t''_0)\cap \tD X,s\Big)\cap\splitset\Big(\pi(t''_1)\cap \tD X,s\Big)\subseteq m. \tag{II$_\laddermap$}
\end{equation*}
\end{lem}

We immediately point out the following observation.
\begin{fact}\label{f.T.extend}
In the context of the lemma above, $\cond$ is equivalent to the following
condition: 
\begin{multline}\label{e.T.extend}
\left( \forall (w_0,w_1) \in  [T^{\tD X,\tD d}_{t'_0}] \times [T^{\tD X,\tD d}_{t'_1}] \right)
\; \\
\splitset\Big(\big(\pi(w_0)\cap \tD X \big)^{++},s\Big)\cap\splitset\Big(\big(\pi(w_1)\cap \tD X \big)^{++},s\Big)\subseteq m.
\end{multline}
\end{fact}
\begin{proof}
This follows mechanically, if somewhat tediously, from the definitions. We give some hints:
For the equivalence of \eqref{i.main} with \eqref{e.T.extend} when $\laddermap(s) > 1$, observe that 
\[
t'_i \in T^{\tD X, \tD d\cup\{s\conc m'\}}
\]
if and only if there is $w_i \in [T_{t'_i}]$ such that
\begin{equation}\label{e.d.u.sm'}
\tD d \cup \{ s \conc m'\} \subseteq \big(\pi(w_i)\cap \tD X\big)^{++}.
\end{equation}
As $s \in \tD d$, Equation~\eqref{e.d.u.sm'} is equivalent to
\[
t'_i \in T^{\tD X, \tD d} \; \land\; m' \in \splitset\Big(\big(\pi(w_1)\cap \tD X \big)^{++},s\Big).
\]
For the equivalence of \eqref{i.main.cond} with \eqref{e.T.extend} in the case that $\laddermap(s) = 1$, i.e., when $s\conc m' \in S_\laddermap$, notice that
\[
\splitset\Big(\big(\pi(w_1)\cap \tD X \big)^{++},s\Big) = \splitset\Big(\pi(w_1)\cap \tD X,s\Big),
\]
since $s \in \big(\pi(w_1)\cap \tD X \big)^{++}$.
Further details are elementary.
\end{proof}

\begin{rem}
For the reader's orientation, we point out that the differences between the lemma above and
its two-dimensional cousin, Lemma~\ref{l.2D.branch} are not great:
This is manifest for \eqref{i.2D.branch.I}  and \eqref{i.main}.
In the context of Item \eqref{i.2D.branch.II} of Lemma~\ref{l.2D.branch},
since
$m_0 \in d$, it holds that
\begin{equation*}
\big(\pi(t''_0)\cap \pi(t''_1) \cap \tD X \big)(m_0) =
\splitset\Big(\pi(t''_0)\cap \tD X,\langle n\rangle\Big) \cap \splitset\Big(\pi(t''_1)\cap \tD X,\langle n\rangle\Big).
\end{equation*}
Moreover, $m_0 \in d$ could be rewritten as $s:= \langle m_0 \rangle \in \tD d$, using the trivial correspondence of $\nat$ with $\nat^1$, to match the notation of this section.
As a consequence, Lemma~\ref{l.2D.branch} is identical (modulo notational differences) with the special case of 
Lemma~\ref{l.alphaD.branch} above where $\gamma$ is the rank function on $U = \nat^{\leq 2}$.
\end{rem}

\medskip

For the proof of the lemma we need a higher dimensional analogue of $\sqsubset_2$.
Recall that since $\higherdom(S_\laddermap) \cup S_\laddermap$ is a tree (a subtree of $\nat^{<\nat}$ ordered by $\subseteq$),
we can speak of a \emph{subtree of $\higherdom(S_\laddermap) \cup S_\laddermap$} 
with the usual, obvious meaning.

\begin{dfn}
Given two finite subtrees $\tD c, \tD c'$ of $\higherdom(S_\laddermap) \cup S_\laddermap$ 
write $\tD c \sqsubset_\laddermap \tD c'$ to mean that
$\tD c \subseteq \tD c'$ and
 for each $s\in \tD c$,
\[
\sup \big( \splitset(\tD c,s) \big) < \sup  \big( \splitset(\tD c',s) \big),
\]
that is,
\[
\{n\in\nat \setdef s \conc n \in \tD c\} \sqsubsetneq \{n\in\nat \setdef s \conc n \in \tD c'\}
\]
(where, as we wish to remind the reader, for $a,a'\in[\nat]^{<\infty}$,  $a \sqsubsetneq a'$
 means that
$\hat a$ is proper initial segment of $\hat a'$).
 \end{dfn} 
 By the following obvious fact (which we state without proof), this gives us a way to build sets in the co-ideal:
 \begin{fact}\label{f.alphaD.grow}
 If $\langle \tD c_k \setdef k \in \nat\rangle$ is a sequence consisting of subtrees of 
 $\higherdom(S_\laddermap) \cup S_\laddermap$ and
$\tD c_k \sqsubset_\laddermap \tD c_{k+1}$ for each $k\in \nat$, 
 it must hold that
 \[
 S_\laddermap \cap \bigcup_{k\in\nat} \tD c_k \in \finalphpp.
 \]
 \end{fact}
 
 The proof of Lemma~\ref{l.alphaD.branch'} above can now be carried out closely following the blueprint of the proof of its two-dimensional analogue, Lemma~\ref{l.2D.branch}.
 \begin{proof}[Proof of Lemma~\ref{l.alphaD.branch'}]
 Again, assume $\tD X \in \finsqp$ and
fix $t_0, t_1 \in T^{\tD X}$ such that $\lh(t_0) = \lh(t_1)$ but $\pi(t_0)\neq \pi(t_1)$.
As before we suppose the lemma fails and derive a contradiction.

Analogously to the 2-dimensional case, We build sequences
\begin{itemize}
\item $\emptyset \neq \tD{d}^0 \subseteq \tD{d}^1 \subseteq  \tD{d}^2 \subseteq \hdots$ from $[\higherdom(S_\laddermap)]^{<\infty}$,
\item $\tD{b}^0  \subseteq  \tD{b}^1  \subseteq  \tD{b}^2  \subseteq  \hdots$ from $[S_\laddermap]^{<\infty}$,
\item $t^0_i  \subseteq  t^1_i  \subseteq  t^2_i  \subseteq \hdots$ from $T_{t_i}$ for each $i\in\{0,1\}$, 
\end{itemize}
such that %
for each $k\in \nat$,
\begin{enumerate}[label=(\Alph*)]
\item\label{IH.grow'}  $\tD{d}^k \cup \tD{b}^k \sqsubset_\laddermap \tD{d}^{k+1}\cup \tD{b}^{k+1}$,
\item\label{IH.subseteq'} $\tD{b}^{k+1} \subseteq \pi(t^k_0)\cap \pi(t^k_1)$,
\item\label{IH.keep.going'} $t^{k+1}_i \in T^{\tD X,\tD{d}^k}_{t^k_i}$.
\end{enumerate}
In the end, we will have $\bigcup_{k \in \nat} \tD{d}^k = \higherdom\left(\bigcup_{k \in \nat} \tD{b}^k\right)$.

Let $\tD{d}^0 = \{\emptyset\}$ (this will ensure that the following construction's first step, i.e., the case  $k=0$, is not vacuous), $\tD{b}^0=\emptyset$ and $t^0_i = t_i$ for each $i\in\{0,1\}$.

We describe the inductive step of the construction:
Suppose we have constructed $\tD{d}^k$, $\tD{b}^k$, $t^k_0$, and $t^k_1$.
To find $\tD{b}^{k+1}$ satisfying \ref{IH.grow} we must repeatedly apply the assumption that Lemma~\ref{l.2D.branch} fails. 
We make finitely many such applications---to be precise, one for each $s \in \tD{d}^k$.

Write $\bar{j} = \lvert \tD{d}^k \rvert$  and let $\langle s^j \setdef j < \bar{j}\rangle$ enumerate  $\tD{d}^k$. 
We build finite sequences 
$\langle \tD{d}^{k,j} \setdef j \leq \bar{j}\rangle$,  $\langle \tD{b}^{k,j} \setdef j \leq \bar{j}\rangle$, 
and $\langle t^{k,j}_i \setdef j \leq \bar{j}\rangle$ for $i \in\{0,1\}$, 
starting with $\tD{d}^{k,0} = \tD{d}^k$, $\tD{b}^{k,0}=\tD{b}^k$, and $t^{k,0}_i = t^k_i$.

Suppose now we already have  $\tD{d}^{k,j}$, $\tD{b}^{k,j}$, and $t^{k,j}_0$ and $t^{k,j}_1$.
By, assumption, the statement of the lemma fails. 
Therefore, taking
$\tD{d} = \tD{d}^{k,j}$, 
$t'_i = t^{k,j}_i$ for each $i \in \{0,1\}$, 
$s = s^j$, and $m$ such that
\[
\splitset (\tD{d}^k \cup \tD{b}^k,s)\subseteq m
\] 
the condition $\cond=\condargs$ fails. We argue by cases.

First suppose that $\laddermap(s) >1$ and \eqref{i.main} fails. 
We can find $m' \geq m$ 
such that
\[
(\forall i \in \{0,1\})\;   t'_i \in T^{\tD X, \tD d \cup \{ s \conc m'\}}.
\]
We can let
$\tD{d}^{k,j+1} = \tD{d}\cup\{s\conc m'\}$ 
and
$t^{k,j+1}_i = t'_i = t^{k,j}_i$;
we have maintained \ref{IH.keep.going'} and are one step closer to
ensuring \ref{IH.grow'}.
Finally, let $\tD{b}^{k,j+1} = \tD{b}^{k,j}$; then \ref{IH.subseteq'} is trivial by induction. 

Now consider the case that $\laddermap(s)=1$ and
\eqref{i.main.cond} fails. 
We can find 
 $m' \geq m$ together with 
 \[
 (t''_0, t''_1) \in  T^{\tD X,\tD d}_{t'_0} \times T^{\tD X,\tD d}_{t'_1}
 \]
  such that
 \[
m' \in  \splitset\Big(\pi(t''_0)\cap \tD X,s\Big)\cap\splitset\Big(\pi(t''_1)\cap \tD X,s\Big).
 \]
In this case let $t^{k,j+1}_i = t''_i$ for each $i\in\{0,1\}$, $\tD{b}^{k,j+1} = \tD{b}^{k,j}\cup\{s \conc m'\}$, and $\tD{d}^{k,j+1} = \tD{d}$, which takes us one step closer to
ensuring \ref{IH.grow'}
without damaging \ref{IH.subseteq'}.
Moreover we maintain  \ref{IH.keep.going'} by choice of $t''_i$.
This finishes the inductive step in the construction of 
$\langle \tD{d}^{k,j} \setdef j \leq \bar{j}\rangle$,  $\langle \tD{b}^{k,j} \setdef j \leq \bar{j}\rangle$, 
and $\langle t^{k,j}_i \setdef j \leq \bar{j}\rangle$ for $i \in\{0,1\}$.

Finally let 
$\tD{d}^{k+1} = \tD{d}^{k,\bar{j}}$, 
$\tD{b}^{k+1} = \tD{b}^{k,\bar{j}}$ and $t^{k+1}_i = t^{k,\bar{j}}_i$ for each $i\in\{0,1\}$.
It is clear from the inductive construction that \ref{IH.subseteq'} and  \ref{IH.keep.going'} hold; 
\ref{IH.grow'} holds by each choice of $m'$ and since $\tD d^{k,j+1} \cup \tD b^{k,j+1} = \tD d^{k,j+1} \cup \tD b^{k,j} \cup \{ s_j \conc m'\}$ in the sub-induction over $j$. 

This finishes the inductive step from $k$ to $k+1$ and thus the definition of 
 $\tD{d}^k$, $\tD{b}^k$ and $t^k_i$ for each $k\in\nat$ and $i\in\{0,1\}$.
But (exactly as in the 2-dimensional case)
$\langle t^k_i \setdef k\in\nat\rangle$ gives us a branch $w_i \in [T_{t_i}]$ for each $i\in\{0,1\}$ such that
$\pi(w_0)\neq\pi(w_1)$ and 
\[
\bigcup_{k\in\nat} \tD{b}^k \subseteq \pi\left(w_0\right)\cap \pi\left(w_1\right)
\]
which contradicts that $\mathcal A$ is a $\finalph$-almost disjoint family since the left hand side is an element of 
$\finalphp$ by \ref{IH.grow'} and Fact~\ref{f.alphaD.grow}. Having reached a contradiction, the proof of Lemma~\ref{l.alphaD.branch'} is complete.
 \end{proof}

\subsection{Proof of the higher dimensional theorem}

We are now ready to state and prove the higher dimensional analogue of Theorem~\ref{t.2D.analytic}.

\begin{thm}\label{t.alphaD.analytic}
Suppose $\mathcal A$ is an analytic $\finalph$-almost disjoint family  and $\vec{z} = \langle \tD{Z}^l\setdef l\in\nat\rangle$ is a sequence of pairwise disjoint sets from $\finalphpp$ satisfying 
\ref{Z.refine.alph} and \ref{Z.disjoint.alph} above. 
For any $W \in [\nat]^{\infty}$ there is $W_0 \in [W]^{\infty}$ such that  $\fromx[\vec{z}]{W_0}$ is $\finalph$-a.d.\  from each element of $\mathcal A$.
\end{thm} 
\begin{proof}
The proof strategy is similar to the 2-dimensional case: 
Let $\mathcal A$ and $\vec{z}$ as in the lemma be given.
Fix a tree $T$ on $2\times\nat$ such that $\pi[T]=\mathcal A$ (modulo our usual identification of $2^\nat$ with
$\powerset(\nat^{<\nat})$).

Similarly to the 2-dimensional case, we show how to find a ``large set'' on which for each $\tD d \in [S_\laddermap]^{<\infty}$ the function 
\begin{gather*}
[\nat]^\infty \to \text{set of trees on $2\times \nat$},\\
B \mapsto T^{\from[B],\tD d}
\end{gather*}
is
constant,
using the fact that analytic sets are completely Ramsey and ``invariance'' of each of the trees.
\begin{claim}\label{c.alphaD.constant}
For every $W \in [\nat]^{\infty}$ there is $W_0 \in [W]^{\infty}$ such that 
\begin{equation}\label{e.W_0}
(\forall W' \in [W_0]^{\infty})(\forall \tD d \in [\higherdom(\from[W'])]^{<\infty})\; T^{\from[W'],\tD d}=T^{\from[W_0],\tD d}
\end{equation}
\end{claim}
\begin{proof}[Proof of Claim~\ref{c.alphaD.constant}]
The proof is by a fusion argument similar to that of Claim~\ref{c.2D.constant}.
Let $\langle(\tD d_k, t_k)\setdef k\in\nat\rangle$ enumerate 
\[
Z:= \{\tD d \in [S_\laddermap]^{<\infty} \setdef\text{ $\tD d$ is a subtree of $S_\laddermap$}\}^{<\infty}\times T
\] 
so that each element of $Z$ occurs infinitely many times in the enumeration.
Let for each $k\in T$
\[
D_k=\{B \in[\nat]^\infty\setdef t_k \in T^{\from[B],\tD d_k}\}
\]
and note that this set is analytic and hence, completely Ramsey.

Now construct a sequence $(b_k,C_k) \in [\nat]^{<\infty}\times[\nat]^{\infty}$ for $k\in \nat$ and obtain $W_0$ verbatim as in 
the proof of Claim~\ref{c.2D.constant}.

To see that $W_0$ is as desired, let $W'\in[W_0]^{\infty}$ and a finite set $\tD d \subseteq \higherdom(\from[W'])$ which is closed under initial segments be given.
To show $T^{\from[W'],\tD d}=T^{\from[W_0],\tD d}$, consider an arbitrary $t \in T$.
Find $k\in\nat$ such that $(\tD d_k,t_k)=(\tD d,t)$ and $\tD d \subseteq \higherdom(\from[b_k])$.
Since $W' \subseteq_{\fin} C_{k+1}$ we may choose $W''$ such that $W'' \mathbin{E_0} W'$ and $W'' \in [b_k, C_{k+1}]^{\infty}$.
Thus by Equation~\eqref{e.make.constant.one.step},
\[
t \in T^{\from[W''],\tD d} \iff t \in T^{\from[C_{k+1}],\tD d}. 
\]
Moreover, applying conditional invariance of this tree as expressed in Fact\ref{f.alphaD.T^z,d} 
\eqref{i.alphaD.T^z,d.invariant.corollary}, since
$W'' \mathbin{E_0} W'$ and $\tD d \subseteq \higherdom(\from[W''])$ %
we have
\[
T^{\from[W''],\tD d} = T^{\from[W'],\tD d}.
\]
and so also
\[
t \in T^{\from[W'],\tD d} \iff t \in T^{\from[C_{k+1}],\tD d}.
\]
The very same argument, again using higher dimensional conditional invariance from Fact\ref{f.alphaD.T^z,d}\eqref{i.alphaD.T^z,d.invariant.corollary}, also yields
\[
t \in T^{\from[W_0],\tD d} \iff t \in T^{\from[C_{k+1}],\tD d},
\]
proving $T^{\from[W_0],\tD d}=T^{\from[W'],\tD d}$, as $t$ was arbitrary.
\renewcommand{\qedsymbol}{{\tiny  Claim~\ref{c.alphaD.constant}~}$\Box$}
\end{proof}
Let us now fix $W_0 \in[\nat]^{\infty}$ satisfying \eqref{e.W_0} and 
for each 
$\tD d \in \big[\higherdom(S_\laddermap)\big]^{<\infty}$,
write $T^{*,\tD d}$ for $T^{\from[W_0],\tD d}$.
Again, $W_0$ is ``sufficiently generic'' to allow us to show that $\from[W_0]$ is disjoint from every set from $\mathcal A$, via the following two claims:
\begin{claim}\label{c.alphaD.not1branch}
$\lvert \pi[T^{*,\emptyset}]\rvert \neq 1$.
\end{claim}
\begin{proof}[Proof of Claim~\ref{c.alphaD.not1branch}]
The reader may refer to Claim~\ref{c.2D.not1branch} and adapt its proof by simply replacing $\finsq$ by $\finalph$, 
and using %
Facts~\ref{f.from.alphaD}\ref{i.alphaD.tilde.ad}. 
\renewcommand{\qedsymbol}{{\tiny  Claim~\ref{c.alphaD.not1branch}~}$\Box$}
\end{proof}
A detailed account is in order for the next and final claim. 
The major difference to Claim~\ref{c.2D.one.branch} is that we must use the 
Pigeonhole Principle from 
Facts~\ref{f.from.alphaD}\ref{i.alphaD.pidgeon} 
 instead of Facts~\ref{f.from} as well as the higher dimensional branch lemma (Lemma~\ref{l.alphaD.branch}) instead of its 2-dimensional counterpart  Lemma~\ref{l.2D.branch}.
\begin{claim}\label{c.alphaD.one.branch}
$\pi[T^{*,\emptyset}]=\emptyset$.
\end{claim}
\begin{proof}[Proof of Claim~\ref{c.alphaD.one.branch}]
Otherwise by the previous claim, $\lvert \pi[T^{*,\emptyset}]\rvert\geq 2$.
Fix $t_0, t_1 \in T^{*,\emptyset}$ with $\lh(t_0) = \lh(t_1)$ and $\pi(t_0)\neq \pi(t_1)$.
Use Lemma~\ref{l.alphaD.branch} to find $\tD d \in \big[\higherdom(S_\laddermap)\big]^{<\infty}$, $s \in \tD d$, $m\in\nat$, 
and 
$t'_i \in  T^{*,\tD d}_{t_i}$ for each $i\in \{0,1\}$  
such that 
$\condargs$ holds---or equivalently by Fact~\ref{f.T.extend}, \eqref{e.T.extend} holds. 

Let
\[
X_i = \bigcup_{w \in [T^{*,\tD d}_{t'_i}]}  
\splitset\Big(\big(\pi(w_0)\cap \tD X \big)^{++},s\Big)
\]
From \eqref{e.T.extend}, 
we infer
\begin{equation*}
X_0\cap X_1 \subseteq m.
\end{equation*}
Use the Pigeonhole Principle from 
Facts~\ref{f.from.alphaD}\ref{i.alphaD.pidgeon}
 to find $W_1 \in [W_0]^{\infty}$ such that 
\begin{equation}\label{e.d}
\tD d \subseteq \higherdom(\from[W_1])
\end{equation}
and 
\[
\splitset(\from[W_1],s) \subseteq_{\fin} X_0\text{ or }\splitset(\from[W_1],s) \subseteq_{\fin} \nat\setminus X_0.
\]
Just as in the 2-dimensional case, $X_0 \cap X_1 \in \fin$ and so we may fix $i\in\{0,1\}$ so that
\begin{equation}\label{e.alphaD.small}
\splitset(\from[W_1],s)\cap X_i \in \fin.
\end{equation} 
Again, we show that this contradicts $t'_i \in T^{*,\tD d}$:

Now $T^{*,\tD d} = T^{\from[W_1],\tD d}$ 
by choice of $W_0$ and by \eqref{e.d}; and so $t'_i  \in T^{\from[W_1],\tD d}$.
By definition of the tree and since $s\in \tD d$, there is $w \in [T_{t'_i}]$ such that 
$s \in \big(\pi(w) \cap \from[W_1]\big)^{++}$ 
and therefore
\begin{equation}\label{e.alphaD.big}
\splitset\Big(\big(\pi(w)\cap \from[W_1] \big)^{++},s\Big) \text{ is infinite}. 
\end{equation} 
Also by definition of the tree, $w$  is a branch through $T^{W_1,\tD d}_{t'_i}$ and hence through $T^{*,\tD d}_{t'_i}$. 
But by choice of $X_i$ this entails 
\[
\splitset\Big(\big(\pi(w)\cap \from[W_1] \big)^{++},s\Big) \subseteq X_i
\]
 which
  is an absurdity by Equations \eqref{e.alphaD.small} and \eqref{e.alphaD.big}. 
We have reached a contradiction, proving the claim. 
\renewcommand{\qedsymbol}{{\tiny  Claim~\ref{c.alphaD.one.branch}~}$\Box$}
\end{proof}

By the last two Claims it has to be the case that 
$[T^{\from[W_0],\emptyset}]=\emptyset$. 
By Facts~\ref{f.alphaD.T^z,d}(\ref{i.alphaD.T^z,d.pruned}) this means that 
$\emptyset \notin T^{\from[W_0],\emptyset}$ and so by Facts~\ref{f.alphaD.T^z,d}(\ref{i.alphaD.T^z,d.non-empty}) 
we conclude that $\mathcal A$ is not maximal.
\renewcommand{\qedsymbol}{{\tiny  Theorem~\ref{t.alphaD.analytic}~}$\Box$}
\end{proof}

\subsection{Corollaries}\label{s.corollaries}

Just as in the two dimensional case, we can now prove our main theorem for any $\finalph$ as a corollary of the more technical Theorem~\ref{t.alphaD.analytic} which was proved in the previous section.
\begin{cor}\label{c.Finalph}
Theorem~\ref{t.Finalph} holds.
\end{cor}
\begin{proof}
The theorem follows from Theorem~\ref{t.alphaD.analytic} precisely as Theorem~\ref{t.2D.precise} followed from Theorem~\ref{t.2D.analytic} (see Section~\ref{proof.t.2D.precise}, page~\pageref{proof.t.2D.precise}).
\end{proof}

A number of results that were previously shown by somewhat different methods follow quickly from Theorem~\ref{t.Finalph}. 
\begin{cor}[\cite{haga-schrittesser-toernquist}]\label{c.alphaD.no.analytic}
There are no analytic infinite $\mathcal I$-\emph{mad} families, for any $\mathcal I \in \mathfrak F$.
\end{cor}
\begin{proof}
This clearly follows from Theorem~\ref{t.Finalph} (Corollary~\ref{c.Finalph}) above taking $\Gamma=\Sigma^1_1$ and using Jankov-von Neumann uniformization (see \cite[18.1]{kechris}) and the fact that $\sigma(\Sigma^1_1)$ sets have the Ramsey property.
\end{proof}
The special case $\mathcal I = \fin$ (i.e., for classical mad families) of the following corollaries was also shown by Neeman and Norwood \cite{neeman-norwood}.
Moreover, both corollaries were shown in full in \cite{haga-schrittesser-toernquist} via a different proof.
For the definition of $\AD^{+}$ the reader can consult, e.g., \cite{larson-adplus}.
\begin{cor}%
Assuming $\AD^{+}$,
there are no infinite $\mathcal I$-\emph{mad} families, for any $\mathcal I \in \mathfrak F$.
\end{cor}
\begin{proof}
By definition  $\AD^{+}$ implies $\mathsf{DC}_\reals$. 
To prove the corollary, by Theorem~\ref{t.Finalph} it suffices to show that $\RUnif^*(\Gamma)$ holds for $\Gamma = \powerset(\nat^\nat)$.
That every set is Ramsey measurable follows from the results in \cite{ikegami},
while $\RUnif(\Gamma)$ follows from the results in \cite{muller}.
Together, these are equivalent to $\RUnif^*(\Gamma)$.\footnote{That $\AD^+ \Rightarrow \RUnif^*(\Gamma)$ may have been known for a long time; but alas, we do not know any earlier reference in the literature.}
\end{proof}

\begin{cor}[\cite{haga-schrittesser-toernquist}]%
Let $\mathcal I \in \mathfrak F$.
Under $\PD + \DCR$, there are no projective $\mathcal I$-\emph{mad} families.
\end{cor}
\begin{proof}
By $\PD + \DCR$ every projective subset of $[\nat]^\infty$ has the Ramsey property, and projective uniformization (see \cite{kechris}) holds, so
$\RUnif^*(\Gamma)$ holds with $\Gamma$ the pointclass of all projective sets. 
Therefore, the corollary follows by Theorem~\ref{t.Finalph}.
\end{proof}

\begin{cor}[\cite{toernquist, haga-schrittesser-toernquist}]%
Suppose there is an inaccessible cardinal.
In Solovay's model, here are no $\mathcal I$-\emph{mad} families, for any $\mathcal I \in \mathfrak F$.
\end{cor}
\begin{proof}
It is well known that $\DC$ as well as $\RUnif^*(\Gamma)$ hold for $\Gamma = \powerset(\nat^\nat)$ in Solovay's model.
The rest of the argument is as in the previous corollary.
\end{proof}

\bibliography{schrittesser-toernquist-ramsey-property-higher-mad}{}
\bibliographystyle{amsplain}
 
\end{document}